%% file: main.tex
\documentclass{article}
\usepackage[a4paper,left=2.cm,right=2.cm,top=3.75cm,bottom=3.75cm]{geometry}

\usepackage{graphicx} 
\usepackage{amsmath,amssymb,amsthm}

\usepackage{subfig}
\usepackage{xcolor}
\usepackage{adjustbox}
\usepackage{url}
\usepackage{mathrsfs}
\usepackage{amsbsy}
\usepackage{multirow}
\usepackage{tikz}

\usetikzlibrary{spy}
\usepackage[most]{tcolorbox}

\newtheorem{theorem}{Theorem}
\newtheorem{corollary}[theorem]{Corollary}

\newtheorem{prop}[theorem]{Proposition}
\newtheorem{definition}{Definition}

\def\N{\mathbb{N}}
\def\R{\mathbb{R}}

\DeclareMathOperator*{\sign}{sign}

\newcommand{\mcR}{\mathcal{R}} 
\newcommand{\mcRG}{\mathcal{R}_{\text{lim}}} 
\newcommand{\mcRS}{\mathcal{R}_{\text{sparse}}} 
\newcommand{\mcRA}{\mathcal{R}_A} 
\newcommand{\mcLG}{\mathcal{L}_{\text{lim}}} 
\newcommand{\mcLS}{\mathcal{L}_{\text{sparse}}} 
\newcommand{\mcLA}{\mathcal{L}_A} 
\newcommand{\mcS}{\mathcal{S}}

\renewcommand{\vec}[1]{\boldsymbol{#1}} 
\newcommand{\mb}{\vec{m}} 
\newcommand{\ub}{\vec{u}} 
\newcommand{\wb}{\vec{w}} 
\newcommand{\zetab}{\vec{\zeta}} 
\newcommand{\RbA}{\vec{\mathcal{R}}_A} 
\newcommand{\W}{\mathcal{W}}
\newcommand{\Wb}{\vec{\mathcal{W}}} 

\newcommand{\PDONetB}{\boldsymbol{\Lambda}_{\zetab_k}} 
\newcommand{\PDONet}{\Lambda_k} 

\title{Revisiting $\Psi$DONet: microlocally inspired filters for incomplete-data tomographic reconstructions}
\author{
    Tatiana A.~Bubba\thanks{Department of Mathematics and Computer Science, University of Ferrara, Italy ({tatiana.bubba@unife.it})} 
    \and Luca Ratti\thanks{Department of Mathematics, University of Bologna, Italy ({luca.ratti5@unibo.it})} 
    \and Andrea Sebastiani\thanks{Department of Physics, Computer Science and Mathematics, University of Modena and Reggio Emilia, Italy ({andrea.sebastiani@unimore.it})} 
}
\date{}

\begin{document}

\maketitle

\begin{abstract}
In this paper, we revisit a supervised learning approach based on unrolling first introduced in~\cite{bubba2021deep} and called $\Psi$DONet, by
providing a deeper microlocal interpretation for its theoretical analysis, and extending its study to the case of sparse-angle tomography. Furthermore, we refine the implementation of the original $\Psi$DONet considering special filters whose structure is specifically inspired by the streak artifact singularities characterizing tomographic reconstructions from incomplete data. This allows to considerably lower the number of (learnable) parameters while preserving (or even slightly improving) the same quality for the reconstructions from limited-angle data and providing a proof-of-concept for the case of sparse-angle tomographic data.
\end{abstract}


\section{Introduction}
In~\cite{bubba2021deep}, the authors introduced $\Psi$DONet, a supervised approach based on a convolutional neural network (CNN) designed starting from the unrolled iterations of the Iterative Soft Thresholding Algorithm (ISTA)~\cite{daubechies2004iterative}, a classical  algorithm applied to the solution of a sparsity promoting regularization formulation for inverse problems (cf.~Equation \eqref{eq:sparsprom}). 
In particular, $\Psi$DONet was applied to the inverse problem of limited-angle tomography, a limited-data tomographic problem where only a subset (or wedge) of the angular views is available (cf.~Figure~\ref{fig:CucumberTomo}), making it a severely ill-posed inverse problem~\cite{davison1983ill}. The novelty of the approach in~\cite{bubba2021deep} is that $\Psi$DONet takes advantage of a sound theoretical backbone based on the convergence analysis of sparsity promoting regularization~\cite{grasmair2008sparse,grasmair2011necessary} and the microlocal properties of the X-ray transform~\cite{quinto1993singularities,quinto2017artifacts}, which provides a clear interpretation of the singularities and artifacts appearing when reconstructing from Computed Tomography (CT) data. 
In particular, in~\cite{bubba2021deep} the convergence of $\Psi$DONet was proved using classical results~\cite{flemming2018injectivity} on the convergence theory of ISTA and relying on interpreting $\Psi$DONet as a modification or perturbation of ISTA (cf.~\cite[Section 4.3]{bubba2021deep}). Moreover, it was suggested that $\Psi$DONet induces a 
``microlocal regularization'' to support the evidence that the network intervenes to mitigate the propagation of the streak artifacts  singularities, rather than inferring the missing information as in other approaches proposed in the literature (cf.~\cite[Section 4.2]{bubba2021deep}).

In the last decade, there has been an increasing effort in applying data-driven approaches to incomplete-data tomographic problems and more specifically to limited-angle tomography. A very recent review article focussing on CT problems, including limited-angle and sparse-angle tomography, is \cite{kiss2024benchmarking}, where it is provided a benchmarking study on a wide range of algorithms representative for different categories of learned reconstruction methods on a dataset of experimental CT measurements. The class of approaches most relevant to our work is that of methods 
explicitly leveraging the theory of microlocal analysis to assist and guide data-driven strategies for limited-angle tomography~\cite{andrade2022deep,andrade2019extraction,arndt2023model,bubba2019learning,rautio2023learning,wang2024robustness}, thus enabling to add a theoretical layer to the development of the numerical method. 
More in general, the rationale underpinning $\Psi$DONet fits into the line of research of hybrid reconstruction frameworks for the solution of imaging inverse problems (see, e.g., ~\cite{arridge2019solving,bubba2024data} and references therein). This paradigm leverages the practical advantages of learning-based method  with the theoretical understanding that comes from model-based approaches, developing methods that are able to significantly surpass both model- and most pure data-based methods. 

In this work, we revisit $\Psi$DONet, providing a deeper microlocal interpretation for its theoretical analysis, and extend its study to the case of sparse-angle tomography, another example of incomplete-data tomography arising from considering a scarce, equidistant sampling for the angular views (cf.~Figure~\ref{fig:CucumberTomo}). In particular, we introduce a new continuous formulation of $\Psi$DONet which, thanks to the very recent analysis provided in~\cite{andrade2022deep}, allows to argue that: (1) $\Psi$DONet \textit{can} allow to recover edges and discontinuities which are invisible from the incomplete sinogram, although we cannot predict their exact location; (2) $\Psi$DONet not only can dampen the strength of  streak artifacts but \textit{can} prevent them from appearing in the reconstruction.
From a numerical point of view, we refine the implementation of the original $\Psi$DONet in~\cite{bubba2021deep} considering special filters whose structure is specifically inspired by the streak artifact singularities characterizing tomographic reconstructions from incomplete data. This allows, in particular, to considerably lower the number of (learnable) parameters while preserving (or even slightly improving) the same quality for the reconstructions from limited-angle data and providing a proof-of-concept for the case of sparse-angle tomographic data. 

The remainder of this paper is organized as follows. In Section~\ref{sec:MathsPrelim} we review some background material on microlocal analysis. Then, in Section~\ref{sec:PsiDONet} we recall some concepts on the sparse regularization theory as a starting point to introduce the $\Psi$DONet architecture. Section~\ref{subsec:PSIDOjust} contains our main theoretical results, and introduces three variants of $\Psi$DONet, specifically tailored to limited- and sparse-angle tomography. Section~\ref{sec:NumericalExperiments} demonstrates the potentiality of our new microlocal-informed filters through numerical experiments on limited-angle and sparse-angle CT reconstruction. Finally, Section~\ref{sec:Conclusions} summarizes our findings and discusses future research directions.

\section{Mathematical Preliminaries}\label{sec:MathsPrelim}
In this section, we summarize some 
concepts and results from microlocal analysis that are essential for our proposed method. In particular, we recall how microlocal analysis can  
explain which edge information is available in the acquired tomographic data (for more details see, e.g., \cite{frikel2013characterization,nguyen2015strong,quinto2017artifacts} and references therein). 

\subsection{Some microlocal analysis concepts}\label{sec:Microlocal}
One of the main goals of microlocal analysis is to give a precise description of the discontinuities of functions and their propagation through the application of certain families of operators~\cite{hormander1983analysis,krishnan2015microlocal,salo2020applications}. Since these insights are key to our proposed architecture we will briefly summarize some of those in the following. First, the notions of singular support and wavefront sets are required. In the following we denote by $\mathcal{F}$ the Fourier transform.

\begin{definition}[singular support]
A function or distribution $u$ is $C^\infty$ \textit{near} $x_0 \in \R^n$ if there exists $\varphi \in C^\infty_c(\R^n)$ such that $\varphi \equiv 1$ in a neighborhood of $x_0$ and $\varphi (u) \in C^\infty(\R^n)$. Then, we define the singular support $\operatorname{sing\ supp}(u)$ of $u$ as:
    \[
    \operatorname{sing\ supp}(u) = \{ x_0 \in \R^n: \text{ $u$ is not $C^\infty$ near $x_0$} \}.
    \]
\end{definition}

Informally, the singular support of $u$ consists of those points $x_0$ such that $u$ is not a smooth function in any neighbourhood of $x_0$. In essence, the singular support allows to describe the location of a singularity of a function $u$. For instance, if $u$ represents an image, then the singularities in some sense determine the ``sharp features'' of the image.

The wavefront set, introduced next, is a more refined notion of a singularity, because it simultaneously describes the location and the direction of such a singularity. This definition is based on a localized correspondence of smoothness and rapid decay in the Fourier domain.

\begin{definition}[wavefront set]\label{def:WavefrontSet}
A distribution $u$ is microlocally $C^\infty$ near $(x_0,\xi_0)$, $x_0,\xi_0 \in \R^n$, $\xi_0\neq 0$ if there exists $\varphi \in C^\infty_c(\R^n)$ such that $\varphi \equiv 1$ in a neighborhood of $x_0$ and there exists $\psi \in C^\infty_c(\R^n \setminus {0})$ such that $\psi \equiv 1$ in a neighborhood of $\xi_0$ such that $\psi(\xi) \mathcal{F}(\varphi u)(\xi) \leq C_N (1 + |\xi|)^{-N}$ for all $N$. Then, we define the wavefront set $\operatorname{WF}(u)$ of $u$ as:
    \[
    \operatorname{WF}(u) = \{ (x_0,\xi_0) \in \R^n \times (\R^n \setminus \{0\}): \text{ $u$ is not microlocally $C^\infty$ near $(x_0,\xi_0)$} \}.
    \]
\end{definition}

For example, if $u = \chi_D$ is the characteristic function of a bounded, strictly convex $C^\infty$ domain $D$ and if $x_0 \in \partial D$, one could think that $u$ is in some sense smooth in tangential directions at $x_0$ (since $u$ restricted to a tangent hyperplane is identically zero, except possibly at $x_0$), but that $u$ is not smooth in normal directions at $x_0$ since in these directions there is a jump. The singular support of $u$ is clearly $\partial D$.

Next, we introduce the definition of a pseudodifferential operator.

\begin{definition}[Pseudodifferential operators]
A pseudodifferential operator ($\Psi$DO) of order $m$ is any operator $\mathcal{A}$ of the form
    \[
    \mathcal{A}u(x) = (2\pi)^{-n} \int_{\R^n} a(x,\xi) e^{ix\cdot\xi} \widehat{u}(\xi) \, \text{d}\xi,
    \]
    where the symbol $a(x,\xi) \in C^{\infty}(\R^n \times \R^n)$ is such that, for all $\alpha,\ \beta \in \N^n$ 
    \[
    |\partial^\alpha_x \partial^\beta_\xi a(x,\xi)| \leq C_{\alpha,\beta} (1 + |\xi|)^{m-|\beta|} \qquad \forall \xi \in \R^n.
    \]
We denote by $\Psi^m$ the set of $\Psi$DOs corresponding to the symbol class $S^m$, that is, for $a \in S^m$ we write $\mathcal{A} \in \Psi^m$.
\end{definition}

A classic example are \textit{elliptic operators}, that is, operators whose principal symbol (i.e. the part containing the highest order derivatives) $\sum_{|\alpha|= m} a_\alpha(x) \xi^\alpha$ does not vanish for $\xi \neq 0$.

A key fact about a $\Psi$DO is that applying it to a function (or distribution) never creates new singularities. 

\begin{prop}[Pseudolocal and microlocal property of $\Psi$DOs]
Any $\mathcal{A} \in \Psi^m$ has the pseudolocal property 
\[
\operatorname{sing\ supp}(\mathcal{A}u)\subset \operatorname{sing \ supp}(u),
\] 
and the microlocal property
\[
\operatorname{WF}(\mathcal{A}u) \subset \operatorname{WF}(u).
\]
\end{prop}

In particular, elliptic operators are those that completely preserve singularities, i.e., if $\mathcal{A} \in \Psi^m$ is elliptic then for any $u$,
\begin{align*}
\operatorname{sing\ supp}(\mathcal{A}u) &= \operatorname{sing \ supp}(u), \\
\operatorname{WF}(\mathcal{A}u) &= \operatorname{WF}(u).
\end{align*}
Thus, any solution $u$ of $\mathcal{A}u = f$ is singular precisely at those points where $f$ is singular.

A more general class of operators (containing $\Psi$DOs) is that of Fourier Integral Operators (FIOs), taking the form:
\begin{equation}\label{eq:defFIOs}
\mathcal{A}u(x) = (2\pi)^{-n} \int_{\R^n} e^{i \varphi(x,\xi)}a(x,\xi)\widehat{u}(\xi) \, \text{d}\xi,
\end{equation}
where $a(x,\xi) \in S^m$, and $\varphi(x,\xi)$ is a real valued phase function. For $\Psi$DOs the phase function is always $\varphi(x,\xi) = x \cdot \xi$, but for FIOs the phase function can be quite general (for a formal definition and more details on FIOs see~\cite{hormander1983analysis}).

An important property of FIOs is that they, unlike $\Psi$DOs, can move singularities. In brief, any FIO has an associated \textit{canonical relation} that describes what the FIO does to singularities. The canonical relation of the FIO $\mathcal{A}$ defined in~\eqref{eq:defFIOs} is 
\[
C = \{ (x, \nabla_x \phi(x,\xi), \nabla_{\xi}\phi(x,\xi),\xi) \ :  \quad (x,\xi) \in T^*\R^n \setminus \{0\} \},
\]
where $T^*\R^n \setminus \{0\} := \{(x, \xi) \ : \ x, \xi \in \R^n, \xi \not= 0\}$. Then, $\mathcal{A}$ moves singularities according to the rule:
\[
\operatorname{WF}(\mathcal{A}u) \subset C(\operatorname{WF}(u)),
\]
where
\[
C(\operatorname{WF}(u)) := \{ (x,\xi) \ : \ (x,\xi,y,\eta) \in C \quad \text{for some} \; (y,\eta) \in \operatorname{WF}(u)\}.
\]

\subsection{Microlocal analysis for tomographic imaging}\label{sec:MicroTomo} 

We now review some of the main results related to the microlocal analysis of the operators involved with the (parallel-beam) Computed Tomography (CT) application, considering in particular the Radon transform in the plane or the X-ray transform in 2D.
Our main reference is \cite{borg2018analyzing}, which actually relies on relevant results developed in the literature, see \cite{frikel2013characterization,katsevich1997local,nguyen2015strong,quinto1993singularities}.

Recall that, given a function $u \in L^1(\R^2)$, 
the X-ray transform of $u$ is given by:
\begin{equation}\label{eq:RadonTransf}
\mcR (u) (s,\omega):= \int_{-\infty}^{\infty} u(s\omega +t \omega^\perp) \text{d}t, \qquad s\in \R, \ \omega \in S^1.
\end{equation}
Here, $\omega^\perp$ is the vector in $S^1$ obtained by rotating $\omega$ counter-clockwise by $90^\circ$. In particular, the X-ray transform $\mcR$ is an elliptic FIO of order $-1/2$, i.e., singularities are moved along specific lines, whereas the normal operator $\mcR^* \mcR$ is a classical $\Psi$DO of order $-1$ in $\R^2$:
\[
\mcR^* \mcR(u) = \int_{\R^2} \frac{1}{|\xi|} e^{i\xi\cdot x} \widehat{u}(\xi) d\xi.
\]
This in particular means that, for the normal operator, the singularities are recovered:
\[ 
\operatorname{sing\ supp}(\mcR^* \mcR u) = \operatorname{sing\ supp}(u).
\]
Because $\mcR^* \mcR$ is a $\Psi$DO of order $-1$ (i.e., smoothing of order 1), $\mcR^* \mcR u$ gives a slightly blurred version of $u$ (backprojection of the data) where, however, the main singularities are still visible. To make singularities sharper, one can use the Filtered Backprojection (FBP) inversion formula which leverages filtering on top of backprojecting the data. More precisely, the FBP operator can be defined as:
\[
\mathcal{L} u  = \mcR^* \Lambda \mcR u, 
\]
where, for a function $g=g(s,\omega)$, $\Lambda g = \mathcal{F}_s^{-1}(|\xi|\mathcal{F}_s g)$, being $\mathcal{F}_s$ the Fourier transform along the variable $s$.

In the following, we will be mainly concerned with imaging situations where complete X-ray data is not available, in particular, limited-angle and sparse-angle tomography. \textit{Limited-angle tomography} arises when the direction vector $\omega$ is restricted within a limited angular range $[-\Gamma, \Gamma]$, with $\Gamma < \pi/2$: 
\[
\mcRG u := \mcR u_{\big\vert[-\Gamma,\Gamma] \times \R}.
\]
\textit{Sparse-angle tomography} occurs every time the sinogram $\mcR u$ is assumed to be known on a finite set of angles $\{\omega_1,\ldots,\omega_N\} \subset [-\pi/2,\pi/2]$. Nevertheless, to comply with the more general theory proposed in \cite{borg2018analyzing}, we consider the available information restricted within a set with a non-empty interior, consisting of a finite number of strips of width $2\eta$:
\[
\mcRS u := \mcR u_{\big\vert[\omega_1-\eta,\omega_1+\eta] \cup \ldots \cup [\omega_N-\eta,\omega_N+\eta] \times \R}.
\]
We use the notation $\mcRA = \chi_{A \times \R} \mcR$, being $A \subset S^1$ to denote both incomplete-data Radon operators, being $A = [-\Gamma,\Gamma]$ in the case of $\mcRG$ and $A = [\omega_1-\eta,\omega_1+\eta] \cup \ldots \cup [\omega_N-\eta,\omega_N+\eta]$ to represent $\mcRS$. Notice that, in contrast to \cite{borg2018analyzing}, we denote with $A$ the subset of directions in $S^1$ and not of the product space $S^1 \times \R$, as we are only interested in the limited-angle and sparse-angle cases, which are related to an angular subsampling of the X-ray transform.

\begin{figure}
\begin{tabular}{@{}c@{\;}c@{\;}c@{\;}c@{}}
(a) ground truth & (b) full angle & (c) sparse-angle & (d) limited-angle  \\[0.25em]
 \includegraphics[width=0.245\textwidth]{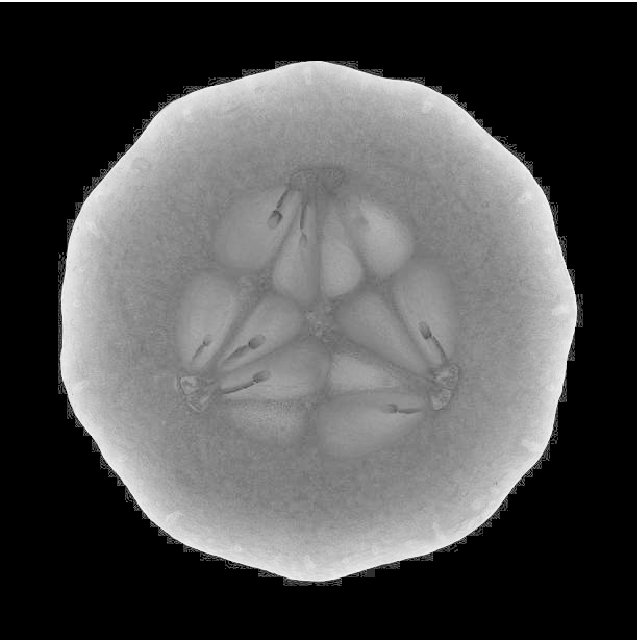}    
    & \includegraphics[width=0.245\textwidth]{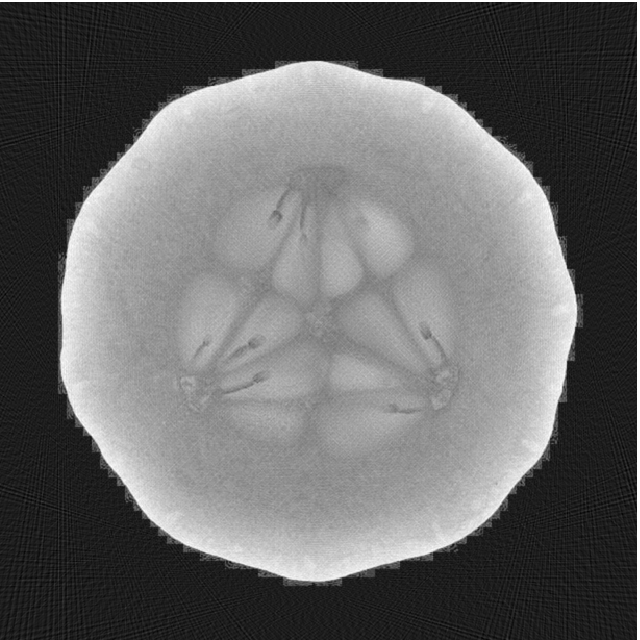} 
    & \includegraphics[width=0.245\textwidth]{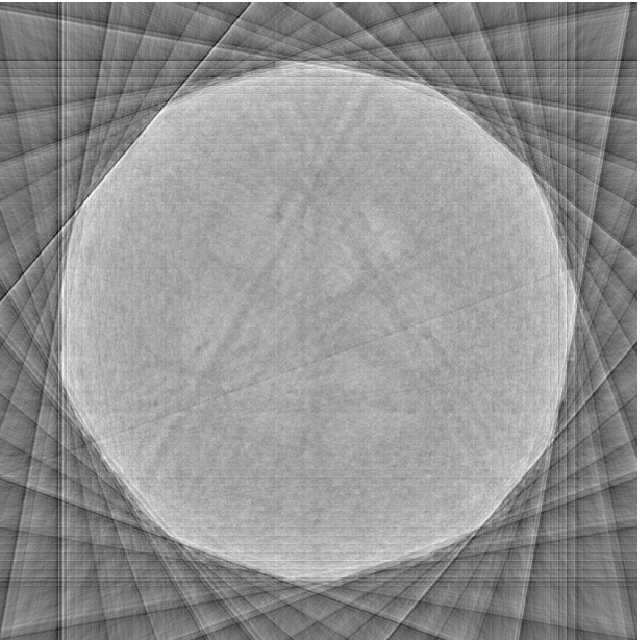} 
    & \includegraphics[width=0.245\textwidth]{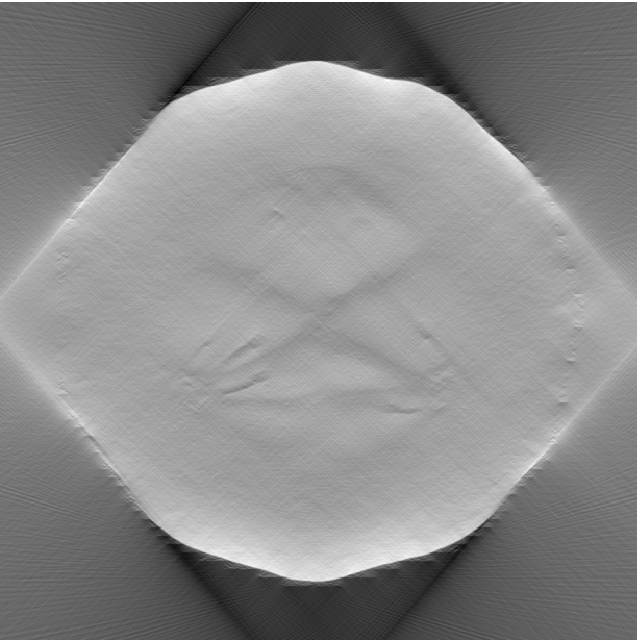} 
\end{tabular}
\caption{Tomographic reconstructions (using FBP) of a cucumber slice from: (b) full data; (c) sparse-angle data; and (d) limited-angle data. Discontinuities along lines are evident when only limited data are are available.}    
\label{fig:CucumberTomo}
\end{figure}

In both cases (limited- and sparse-angle), the normal operator $\mcRA^* \mcRA$ is no longer a $\Psi$DO nor a FIO, but it belongs to a slightly more general class of operators than FIOs~\cite{greenleaf1989non,greenleaf1990estimates}. 
It is nevertheless still possible to have some relevant insights on the propagation of the singularities when $\mcRA^* \mcRA$ is applied. 
In particular, Figure \ref{fig:CucumberTomo} compares the application of $\mathcal{L}$ to an image $u$ with the results obtained by $\mcLG$ and $\mcLS$, where $\mcLA$ is given by:
\[
\mcLA = \mcRA^* \Lambda \mcRA = \mcR^* \Lambda \chi_{A \times \R} \mcR.
\] 
From Figure \ref{fig:CucumberTomo} it is evident that 
the discontinuities of $u$ are perturbed by the application of $\mcLA$ in two ways:
\begin{enumerate}
    \item On the one hand, some portions of the edges, associated with specific normal directions, are not reconstructed.
    \item On the other hand, some additional edges appear, along some specific lines (often referred to as \textit{streak artifacts}).
\end{enumerate}

To precisely describe these effects, one can use microlocal analysis as a powerful tool for predicting which singularities can be recovered from the measurements. In fact, the description of positions and directions of singularities in the sense of Definition~\ref{def:WavefrontSet} allows for a characterization of their visibility in tomographic problems, according to the result below. 

\begin{theorem}[\cite{quinto2017artifacts,quinto1993singularities}]
Let $u \in L^2(\R^2)$ and consider a singularity of $u$, $(x_0,\xi_0) \in \operatorname{WF}(u)$. Take a point $(\omega_0,s_0) \in S^1 \times \R$ such that $x_0 \in L_0$, being $L_0$ the line with direction normal to $\omega_0$ and offset $s_0$, i.e., $L_0 = L(\omega_0, s_0) = \{x \in \R^2: x = s_0 \omega_0 + t \omega_0^\perp, t \in \R \}$. Then,
\begin{itemize}
    \item[(i)]  If $\xi_0$ is normal to $L_0$ ($\xi_0 = \alpha \omega_0$ for some $\alpha \in \R$), then the singularity of $u$ at $x_0$ causes a unique singularity in $\mcR u$ at $(\omega_0, s_0)$.
    \item[(ii)] If $\xi_0$ is not normal to $L_0$, it does not cause any singularity in $\mcR u$ at $(\omega_0, s_0)$. 
\end{itemize}
\label{thm:ArtifactsCharacterization}
\end{theorem}

This can be interpreted by saying that, from the knowledge of the sinogram $\mcR u$ near a point $(\omega_0,s_0)$ it is possible to detect the singularities of $u$ perpendicular to the line $L(\omega_0,s_0)$, but not in other directions. This has a clear implication in incomplete-data scenarios (i.e., where we assume to have access to $\mcRA$, with $A$ being either $A = [-\Gamma,\Gamma]$ or $A = [\omega_1-\eta,\omega_1+\eta] \cup \cdots \cup [\omega_N-\eta, \omega_N + \eta]$), as summarized in the following principle:

\begin{tcolorbox}[colback=gray!5!white,colframe=gray!50!black,
  colbacktitle=gray!15!gray]
\textbf{Principle 1: Visible and invisible singularities}  
\begin{itemize}
    \item If $(x_0, \xi_0)$ is such that $\xi_0 = \alpha \omega_0$, being $\alpha \in \R$ and $\omega_0 \in \operatorname{int}(A)$, then from the singularities of $\mcRA u$ one can determine whether or not $(x_0, \xi_0) \in \operatorname{WF}(u)$. 
    All the singularities of $u$ normal to $L(\omega_0,s_0)$, for $\omega_0 \in A$, $s_0 \in \R$ are called \textit{visible}, as they can be determined from the incomplete sinogram $\mcRA u$.  
    \item
    If $(x_0, \xi_0)$ is such that $\xi_0 = \alpha \omega_0$, but $\omega_0 \notin A$, then the singularities of $\mcRA u$ are not influenced by $(x_0, \xi_0)$.  Such singularities of $u$ (i.e., in directions not contained in $A$), are called \textit{invisible}, and cannot be retrieved from the incomplete sinogram $\mcRA u$.
\end{itemize}
\end{tcolorbox}
An additional consequence is that visible singularities of $u$ appear also in $\operatorname{WF}(\mcLA u)$, whereas invisible ones do not (see \cite[Theorem 3.4]{borg2018analyzing}).
Naturally, which boundaries are \mbox{(in-)visible} is completely determined by the measurement setup, i.e., by the subset $A$. Therefore, it is known a priori and can be used for reconstruction purposes.

To understand the presence of streak artifacts in the reconstruction, the reader is referred to \cite[Theorem 3.7]{borg2018analyzing} and to \cite[Theorem 4.1]{frikel2013characterization} (see also \cite{katsevich1997local,nguyen2015strong}). In particular, in the scenarios we are interested in, the main take-away message is summarized in the principle below.
\begin{tcolorbox}[colback=gray!5!white,colframe=gray!50!black,
  colbacktitle=gray!15!gray]
\textbf{Principle 2: When streak artifacts appear} \\[0.5em]  
    If $\omega_0$ is an angle belonging to $\partial A$, the boundary of the set $A$ (for example, $\omega_0 = \pm \Gamma$ in limited-angle tomography, or $\omega_0 = \omega_1 \pm \eta, \ldots, \omega_N \pm \eta$ in sparse-angle tomography), and $x_0$ is such that $(x_0,\alpha \omega_0) \in \operatorname{WF}(u)$, then any point $(x,\alpha \omega_0)$ with $x \in L(\omega_0,x_0\cdot \omega_0)$ might belong to $\operatorname{WF}(\mcLA u)$.
\end{tcolorbox}
Notice that, as explained in \cite[Section 8.2 and Remark A.5]{borg2018analyzing}, the above result holds true also when the operator $\Lambda$ in $\mcLA$ is replaced by any elliptic pseudodifferential operator, hence is also valid for the normal operator $\mcRA^*\mcRA$.

\section{A pseudodifferential correction via $\Psi$DONet}\label{sec:PsiDONet}
In this section, we recall some key facts on sparse regularization theory for inverse problems, which serves as main motivation to 
introduce the $\Psi$DONet architecture~\cite{bubba2021deep}. 

\subsection{Sparsity-promoting regularization in tomography } 

As already observed in Figure \ref{fig:CucumberTomo}, classical reconstruction methods, such as FBP, do not perform well in incomplete-data scenarios. As an additional complication, in real-world applications the measurements are generally corrupted by noise, i.e., the formulation of the inverse problem associated with $\mcRA$ reads as:
\begin{equation}\label{eq:InvProp}
m = \mcRA u + \epsilon,
\end{equation}
where $\epsilon \in L^2(A\times \R)$ denotes the noise. This induces a perturbation in the reconstruction provided by $\mcLA$ which is extremely amplified: this effect is known as instability. As already observed in \cite{quinto1993singularities}, the invisibility of part of the wavefront set of $u$ is a key factor in describing the stability of its reconstruction.
To compensate this effect and stabilize the inversion, it is common practice to incorporate additional prior knowledge about
the solution into the reconstruction procedure~\cite{benning2018modern}. Specifically, we consider the following variational regularization formulation:
\begin{equation} \label{eq:sparsprom}
\min_u \frac{1}{2} \| \mcRA u - m \|_{L^2}^2 + \lambda \| \W u \|_{\ell^1},
\end{equation}
where $\lambda>0$ is the regularization parameter and $\W: L^2(\R^2) \rightarrow \ell^2$ is a discrete wavelet transform in 2D. We quickly recall that, starting from a suitable mother wavelet $\psi: \R \rightarrow \R$ and its associated scaling function $\phi: \R \rightarrow \R$, it is possible to construct an orthonormal basis of  $L^2(\R^2)$ consisting of wavelets in 2D by rescaling a translating 2D functions obtained by tensor products, as follows:
\[
\left\{\psi_{j,\vec{k}}^{(t)}(x) = 2^j
\psi^{(t)}(
2^j
x - \vec{k}),\quad j \in \mathbb{Z},\ \vec{k} \in \mathbb{Z}^2,\ (t) = (v),(h),(d) \right\},
\]
being 
\[
\begin{gathered} \psi^{(v)}(x_1,x_2) = \phi(x_1)\psi(x_2), \qquad \psi^{(h)}(x_1,x_2) = \psi(x_1)\phi(x_2), \\ \psi^{(d)}(x_1,x_2) = \psi(x_1)\psi(x_2), \qquad \psi^{(f)}(x_1,x_2) = \phi(x_1)\phi(x_2).
\end{gathered}
\]
For all the details, we refer to \cite{mallat1999wavelet}. The resulting family $\{\psi_I\}_I$ (where the index $I$ depends on $j,\vec{k},$ and $t$) creates a discrete basis of the space $L^2(\R^2)$.
In addition, when considering functions on $\Omega \subset \R^2$ (e.g., $\Omega = [0,1]^2$), we can assume that $j\leq J$ (being $J$ the index of the maximum scale, related to the size of $\Omega$) and, for each $j$, we assume a finite set of offsets $\vec{k}$. Finally, if a minimum scale $J_0$ is also fixed (implying that details of size smaller than $2^{-J_0}$ are negligible for applications), we can limit the minimization problem \eqref{eq:sparsprom} to a multiresolution space $W^{J_0} \subset L^2(\Omega)$, namely, the one admitting the following finite-dimensional (orthogonal) basis: 
\[
\Big\{ \psi_{j,\vec{k}}^{(t)}, \quad J_0 \leq j \leq J,\ \vec{k} \in \{1,\ldots,2^j\}^2,\ (t) = (h),(v),(d) \Big\} \cup \Big\{ \psi_{J_0,\vec{k}}^{(f)}, \quad \vec{k} \in \{1,\ldots,2^{J_0}\}^2 \Big\}.
\]
For a fixed basis $\{\psi_I\}_I$, the wavelet transform is such that $[\W u]_I = \langle u ,\psi_I \rangle_{L^2}$. The effect of solving the minimization problem \eqref{eq:sparsprom} is that the solution is expected to provide a good fit with the measured data $m$ and to have only a few non-vanishing components in the wavelet basis, i.e., to have a sparse wavelet representation~\cite{starck2010sparse,alberti2023compressed}. Such a property is quite common in images that are expected to be piecewise constant or smooth with discontinuities in isolated points or lines. Indeed, it is possible to show (see, e.g., \cite{holschneider1995wavelets}) that, in its continuous counterpart, the wavelet transform can resolve the singular support of an image: if no singularity appears at a point, the wavelet coefficients associated with that location decay very rapidly as the scale decreases. Hence, the wavelet coefficient of an image with a few discontinuities is likely to be sparse, which motivates the use of \eqref{eq:sparsprom}. 

Notice that the 2D wavelet transform is not the best tool to encode the sparsity of images whose singularities are concentrated along curves: indeed, no directional information can be encoded with the 2D wavelet transform.
As alternatives, the curvelet and the shearlet transforms 
show the ability to resolve the wavefront set of a distribution (see \cite{candes2005continuous,kutyniok2009resolution}), thus resulting in a sparser representation.
In the context of sparsity-promoting regularization with respect to a curvelet dictionary, we also mention that in   \cite{frikel2013sparse} a precise connection is established between the reconstructed coefficients and the visible singularities of an object, in the context of limited-angle tomography. 
For our purposes, we choose the wavelet transform to leverage a (not optimal, but significant) sparsifying effect while easing the theoretical handling and its associated numerical burden.

The most popular strategy to solve \eqref{eq:sparsprom} is  the Iterative Soft Thresholding Algorithm (ISTA), proposed in \cite{daubechies2004iterative}. Starting from an initial guess $u^{(0)}$, it iteratively creates the subsequent iterates as follows:
\begin{equation}\label{eq:ISTAcont}
    u^{(k+1)} \, = \, \W^* \mcS_{\alpha \lambda}\left(\W (u^{(k)} - \alpha \mcRA^* \mcRA u^{(k)} + \alpha \mcRA^* m) \right),
\end{equation}
where $\alpha < \frac{1}{\| \mcRA^*\mcRA\|}$ and 
$\mcS_\gamma \colon \ell^2 \rightarrow \ell^2$ is the (component-wise) soft-thresholding operator, i.e., $[\mcS_\gamma (w)]_I = S_\gamma(w_I)$, being $ S_\gamma \colon \R \rightarrow \R$ s.t. $S_\gamma(x)= \sign(x) \max(|x| - \gamma, 0)$. The iterates of ISTA are proved to converge (strongly) to a minimizer of \eqref{eq:sparsprom}, see \cite{daubechies2004iterative}. 

It is hard to describe the microlocal properties of a cluster point of \eqref{eq:ISTAcont}. Nevertheless, we observe that the soft-thresholding operator $\mcS_\gamma$ in~\eqref{eq:ISTAcont} might introduce additional singularities in the reconstruction (see also the discussion in Section \ref{subsec:PSIDOjust}). This can be beneficial for the reconstruction, and it might allow to recover part of the invisible singularities. On the other hand, despite the sparsity induced by the regularization term $\| \W u\|_{\ell^1}$, the presence of the operator $\mcRA^*\mcRA$ in Equation \eqref{eq:ISTAcont} might still concur in the formation of streak artifacts, as the numerical evidence suggests.

\subsection{$\Psi$DONet in a nutshell}

Starting from the ISTA iteration~\eqref{eq:ISTAcont}, in~\cite{bubba2021deep} the authors introduced $\Psi$DONet, a convolutional neural network (CNN)
designed to induce a ``pseudodifferential microlocal correction'' that smooths away  singularities like those arising in limited-data tomographic problems, effectively improving the final reconstruction from such limited data. 
In this section, we briefly recall the structure of $\Psi$DONet and what is the rationale behind it.

To begin with, $\Psi$DONet is designed to work on discrete images with a fixed resolution: this requires first to work out a discrete counterpart of Equation~\eqref{eq:ISTAcont}.
To this end, we represent an image $u \in L^2(\R^2)$ by a vector $\ub \in \R^n$ (e.g., a vectorized pixel map), the (incomplete) sinogram as a vector $\mb \in \R^m$, and the limited-data X-ray transform as an operator $\RbA\colon \R^n \rightarrow \R^m$ (or a matrix). The wavelet basis is also truncated at the maximum and minimum scales $J$ and $J_0$, respectively, reflecting the size of the domain and the pixel size, which also implies that the number of wavelet subbands $S$ is fixed. For example, in the case $J-J_0=3$, this implies a total number of $3$ subbands of type $(h),(v),(d)$, and one of type $(f)$, giving $S = 10$. The resulting operator $\Wb: \R^n \rightarrow \R^n$ can be seen as an orthogonal matrix. Next, denoted by $w^{(k)} = \W u^{(k)}$, thanks to the orthogonality of the wavelet operator, the discrete counterpart of \eqref{eq:ISTAcont} can be rephrased as
follows. Starting from an initial guess $\wb^{(0)} = \Wb \ub^{(0)} \in \R^n$, the sequence $\{\wb^{(k)}\}_{k=1}^K$ is generated via the update:
\begin{equation}\label{eq:ISTAitintro}
    \wb^{(k+1)} \, = \, \vec{\mcS}_{\alpha \lambda}\left(\wb^{(k)} - \alpha \Wb \RbA^\top \RbA \Wb^\top \wb^{(k)} + \alpha \Wb \RbA^\top \mb \right),  
\end{equation}
where $\vec{\mcS}_\gamma$ denotes the componentwise soft-thresholding in $\R^n$, i.e., 
$[\vec{\mcS}_\gamma(\wb)]_i = S_\gamma(\wb_i)$.
The solution is eventually recovered by $\ub^{(K)} = \Wb^\top \wb^{(K)}$.
From~\eqref{eq:ISTAitintro}, $\Psi$DONet is obtained  by adding to $\Wb \RbA^\top \RbA \Wb^\top$ a learnable term in the form of a CNN acting as an operator correction. Namely, the unrolled iterations of $\Psi$DONet read as:
\begin{equation}\label{eq:PsiDONet}   
 \wb^{(k+1)} =  \vec{\mcS}_{\gamma_k} (\wb^{(k)} - \alpha_k  \Wb \RbA^\top \RbA \Wb^\top \wb^{(k)} + \alpha_k  \Wb \RbA^\top \mb - \beta_k \PDONetB \wb^{(k)}), 
\end{equation}
where the ``perturbation'' $\PDONetB$ denotes the learnable correction at each layer: $\zetab_k$ are the learnable parameters (i.e., the learnable components of the convolutional filters) of the $k^{\text{th}}$ layer. In addition to the $\zetab_k$'s, the parameters $\alpha_k$, $\beta_k$ and $\gamma_k$ are also learnable.

The novelty underpinning $\Psi$DONet is its structure. $\Psi$DONet is implemented as a combination of upscaling, downscaling and convolution operations, as it is common practice in deep learning, however the way these operations are intertwined and, in particular, the size of its filters is 
inspired by the FIO nature of the operator $\mcR$. With a slight generalization, $\Psi$DONet is tailored to inverse problems whose forward model $\mathcal{A}$ has a normal operator $\mathcal{A}^* \mathcal{A}$ which is a convolutional operator. This is leveraged to provide a convolutional expression of the matrix representing $\mathcal{A}^* \mathcal{A}$ in the wavelet domain (cf.~\cite[Proposition 3.1]{bubba2021deep}), by analytically prescribing the order of convolution, upscaling and downscaling and the filters design, which allows to interpret the operations in~\eqref{eq:PsiDONet} as a layer of a special CNN.   
A possible interpretation of $\Psi$DONet is that the convolutional kernel is split into $\kappa =\kappa_0+\kappa_1$ where $\kappa_0$ is the known part of the model (i.e., $\Wb \RbA^\top \RbA \Wb^\top$) and $\kappa_1$ is an unknown operator (i.e., $\PDONetB$) to be learned, which turns out to be a suitable approximation of a $\Psi$DO. 
As a result, in $\kappa_1$ lays the potential to add information in the reconstruction process with respect to the known part of the model $\kappa_0$, and its $\Psi$DO nature induces a ``microlocal regularization'' that smooths away the singularities (cf.~\cite[Section 4.2]{bubba2021deep}). 
Even though $\Psi$DONet is a CNN designed to reproduce exactly the operations associated with $\Wb \RbA^\top \RbA \Wb^\top$, it provides a fairly general network architecture that allows to recover standard ISTA for a specific choice of the parameters involved. The convergence of the proposed scheme is ensured by classical results on the convergence theory of ISTA and relies on interpreting $\Psi$DONet as a modification or perturbation of ISTA for which the coefficients of the filters can be bounded (cf.~\cite[Section 3.3]{bubba2021deep}).

In the next section, we revisit the convolutional description of $\Psi$DONet, proposing a new continuos formulations that allows for a deeper microlocal interpretation of $\Psi$DONet properties. Then, we introduce three variants of $\Psi$DONet, which use special crafted filters whose structure is inspired by the streak artifact singularities characterizing tomographic reconstructions from incomplete data, including sparse-angle tomography.

\section{$\Psi$DONet: revisiting theory and practice}
\label{subsec:PSIDOjust}

In this section, we return to the theoretical justification provided in \cite{bubba2021deep}, and propose a new
theoretical description of the action of $\Psi$DONet from a microlocal point of view, particularly concerning the propagation of singularities induced by its application. 
To this end, we leverage the approach and some of the results very recently developed in \cite{andrade2022deep}. Compared to the analysis reported in \cite{andrade2022deep}, we limit ourselves to describing the evolution of the singular support of the image, rather than its wavefront set: indeed, since we are considering wavelet-based regularization, the directional information is not available.

From a microlocal point of view, we want to be able to motivate two beneficial effects of the applying $\Psi$DONet: (1) $\Psi$DONet \textit{can} allow to recover edges and discontinuities which are invisible from the incomplete sinogram; (2) $\Psi$DONet \textit{can} prevent streak artifacts from appearing in the reconstruction. 
These results are a significant extension of the microlocal interpretation initially provided in~\cite[Section 4.2]{bubba2019learning}, where it was argued that streak artifacts can be damped thanks to the action of $\Psi$DONet.

\subsection{A novel continuous formulation of $\Psi$DONet}\label{ssec:ContFormulation}
In order to provide a microlocal interpretation of \eqref{eq:PsiDONet}, and to be able to use the tools of microlocal analysis introduced in Section \ref{sec:MathsPrelim}, we should first interpret the $\Psi$DONet algorithm as the discrete formulation of a continuous operator. In  \cite{bubba2021deep}, the $\Psi$DONet algorithm was, in fact, already deduced from a continuous formulation: however, here we provide a different one, which is more convenient for our purposes.

The starting point is to rely on a semi-continuous version of the wavelet transform. In particular, we fix a finite set of scales $j = J_0,\ldots, J$, and, for each subband (identified by $j$ and $(t)$) we consider a continuous offset parameter $z \in \R^2$. This yields a (continuous) translation-invariant dictionary of the form:
\[
\Big\{ \psi_{j,z}^{(t)}(x) =  2^j
\psi^{(t)}( 2^j
(x - z)), \quad J_0 \leq j \leq J,\ z \in \R^2,\ (t) = (h),(v),(d) \Big\} \cup \Big\{ \psi_{J_0,z}^{(f)}, \quad z \in \Omega \Big\},
\]
which results in
a frame, associated with frame bounds and a reconstruction operator, thanks to the theory of Paley-Littlewood transforms (see \cite[Section 5.1.5]{mallat1999wavelet}). Then, the semi-discrete wavelet transform $\mathscr{W}: L^2(\Omega) \rightarrow L^2(\Omega)^Q$ is  defined 
as follows: 
\[
[\mathscr{W}u]_\iota(z) = \int_{\R^2} u(x) \psi_{\iota,z}(x) dx = u \ast \psi_{\iota,z},
\]
where the index $\iota$ denotes the subband (i.e., it depends on $j$ and $(t)$) and $Q$ is the total number of wavelet subbands. The 
adjoint operator $\mathscr{W}^*: L^2(\Omega)^Q \rightarrow L^2(\Omega)$ is defined as:
\[
(\mathscr{W}^* w)(x) = \sum_{\iota \leq Q} \int_{\R^2} w_\iota(z) \psi_{\iota,z}(x)dz .
\]
It is now straightforward to interpret the $\Psi$DONet scheme \eqref{eq:PsiDONet} as the discrete counterpart of the following unrolled iterations:
\begin{equation}\label{eq:PDONetcont}
 w^{(k+1)} \, = \, \pmb{\mathscr{S}}_{\gamma_k}\left(w^{(k)} - \alpha_k \mathscr{W}\mcRA^* \mcRA \mathscr{W}^*w^{(k)} + \alpha_k \mathscr{W}\mcRA^* m - \beta_k \PDONet w^{(k)} \right),
\end{equation}
being $\pmb{\mathscr{S}}_\gamma \colon L^2(\R^2)^Q \rightarrow L^2(\R^2)^Q$ the operator that applies to each subband the operator $\mathscr{S}_\gamma$, i.e., $[\pmb{\mathscr{S}}_\gamma w]_\iota = \mathscr{S}_\gamma(w_\iota)$, where $\mathscr{S}_\gamma \colon L^2(\R^2) \rightarrow L^2(\R^2)$ denotes the Nemytskii operator generated by $S_\gamma \colon \R \rightarrow \R$, namely, $\mathscr{S}_\gamma(u)(x) = S_\gamma(u(x))$. 
Assuming that the translation invariant frame associated with
$\mathscr{W}$ is a Parseval one (which requires some specific assumptions on the mother wavelet $\psi$, see \cite[Theorem 5.11]{mallat1999wavelet}), the synthesis operator $\mathscr{W}^*$ acts also as a reconstruction operator, and we can consider 
$w^{(k)} = \mathscr{W}u^{(k)}$ and $u^{(k)} = \mathscr{W}^*w^{(k)}$.

It is now evident that the singular support of $u^{(k+1)}$ might vary along the iterations 
due to the information coming from the sinogram $m$ and to the application of two operators, $\pmb{\mathscr{S}}_{\gamma_k}$ and $\alpha_k \mathscr{W}\mcRA^*\mcRA\mathscr{W}^* + \beta_k \PDONet$, which we study in more detail in the next section.

\subsection{On introducing edges via soft-thresholding}
Let's consider first
the pointwise soft-thresholding operator $\mathscr{S}_{\gamma}$, applied on each (continuously parametrized) wavelet subband. Our first claim is that $\mathscr{S}_{\gamma}$ might introduce additional singularities in $w^{(k+1)}$, hence in $u^{(k+1)}$. 
To show this, we resort to the techniques introduced in \cite{andrade2022deep} to describe the canonical relation of the Nemytskii operator associated with $\textit{ReLU}(x) = \max\{0,x\}$, based on the following trivial observation:
\begin{equation}
    \label{eq:trivial}
S_{\gamma}(x) = \textit{ReLU}(x-\gamma) - \textit{ReLU}(-x-\gamma).
\end{equation}
We first point out that, in order to interpret the Nemytskii operator associated with $\textit{ReLU}$ as an operator acting on distributions, in \cite{andrade2022deep} a local dampening effect is introduced via a smooth, compactly supported kernel. To ease the notation, here we denote by $\operatorname{ReLU}$ the dampened version 
(which is instead denoted by $\operatorname{ReLU}_{\kappa, \phi_K}$ in \cite{andrade2022deep}), and assume to perform the same operation on $\mathscr{S}_\gamma$, so to interpret it as an operator acting on distributions. Next, we resort to \cite[Theorem 4.11]{andrade2022deep}, which characterizes
the wavefront set of the function $\operatorname{ReLU}(f)$, but we only report here  
the information regarding the singular support, since we are not interested in characterizing the direction of the singularities. This 
significantly simplifies the statement of the original result. In what follows, we denote by $\mathcal{S}$ the Schwartz space of rapidly decreasing functions, by $\mathcal{S}'$ its topological dual, namely, the space of tempered distributions, and by $\langle \cdot,\cdot \rangle_{\mathcal{S}',\mathcal{S}}$ their duality pairing.
\begin{prop}[$\operatorname{sing\ supp}$ version of \protect{\cite[Theorem 4.11]{andrade2022deep}}] \label{prop:andrade}
    Let $\Omega \subset \R^2$ be open and let $u \in \mathcal{S}'(\Omega)$. Then
    \[
    \operatorname{sing\ supp}(\operatorname{ReLU}(u)) \subset \big(\operatorname{sing\ supp}(u) \cap \operatorname{int}(\operatorname{supp}_+(u)) \big) \cup \partial (\operatorname{supp}_+(u)).
    \]
\end{prop} 
Here, $\operatorname{supp}_+(u)$ denotes the \textit{positive support} of a distribution, defined as the closure of the union of all subsets $U$ of $\Omega$ such that $\langle u, \phi\rangle_{\mathcal{S}',\mathcal{S}}>0$ for all $\phi \in \mathcal{S}(\Omega)$ satisfying $\phi \geq 0$ and $\operatorname{supp}\phi \subset U$, excluded $\phi\equiv 0$. In particular, if $u$ is a continuous function, $\operatorname{supp}_+(u) = \operatorname{cl}(\{x: u(x)\geq 0\})$.

An immediate consequence of Proposition \ref{prop:andrade} when applied to~\eqref{eq:trivial} is the following result. Its proof also relies on \cite[Corollary 4.15]{andrade2022deep}, which 
establishes the wavefront set of the sum of two functionals.
\begin{corollary}\label{cor:softthresh}
    Let $\Omega \subset \R^2$ be open and let $u \in \mathcal{S}'(\Omega)$. Then
    \[
    \begin{aligned}
            \operatorname{sing\ supp}(\mathscr{S}_\gamma(u)) \subset \Big ( \operatorname{sing\ supp}(u) \cap \big( \operatorname{int}(\operatorname{supp}_+(u-\gamma)) \cup \operatorname{int}(\operatorname{supp}_+(-u-\gamma))\big) \Big) \\
            \cup\ \partial (\operatorname{supp}_+(u-\gamma)) \cup\ \partial (\operatorname{supp}_+(-u-\gamma)).
    \end{aligned}
    \]
\end{corollary}

The interpretation of the previous result is rather straightforward: the pointwise application of the soft-thresholding operator in \eqref{eq:PDONetcont} may introduce, in each wavelet subband, some additional discontinuities, associated with points where the wavelet coefficients cross the threshold $\pm \gamma$, resulting in some new edges in $u^{(k+1)}$. This is a desirable
consequence, as it might allow recovering some invisible discontinuities, which cannot be reconstructed via $\mcRA^*\mcRA$, although we cannot foresee the location of newly introduced singularities through the proposed analysis.

\subsection{On preventing streak artifacts via kernel smoothing}
\label{ssubsec:kernels}
As discussed in Section \ref{sec:MicroTomo}, the application of the normal operator $\mcRA^*\mcRA$ may cause the insurgence of streak artifacts. In order to avoid this shortcoming, several strategies are possible. For example, in \cite{borg2018analyzing}, 
the authors propose to replace the sharp truncation of the sinogram imposed by $\chi_{A \times \R}$ in $\mcRA$ by a smoothed version, obtained via a function $\bar{\psi} \colon S^1 \times \R \rightarrow \R$ which is equal to $0$ outside $A \times \R$, equal to $1$ in most of $\operatorname{int}(A)\times \R$ and smoothly transitions to $0$ approaching $\partial A\times \R$. The resulting operator $\mcR_{\bar{\psi}} = \bar{\psi} \mcR$ is such that its normal operator is a $\Psi$DO, thus avoiding the generation of additional discontinuities in the reconstruction (cf.~\cite[Theorem 6.1]{borg2018analyzing}).

Also with $\Psi$DONet the idea is to ``correct'' the normal operator $\mcRA^*\mcRA$ is a way that the resulting operator is a $\Psi$DO. However,
the ``microlocal regularizing effect'' powered by $\Psi$DONet is related to a smoothing procedure applied to the Schwartz kernel of $\mcRA^*\mcRA$, represented in the wavelet domain. 
We now formally show this (second) claim.
Recall that the normal operator $\mcR^* \mcR$ of the X-ray transform $\mcR$ is a convolutional operator associated with a Schwartz convolutional  kernel in the Calder\'{o}n-Zygmund family $K(x,y) = \frac{1}{\|x-y\|}= \kappa(x-y)$ for $x \not= y$, whereas it is easy to show that $\mcRA^* \mcRA$ is a convolutional operator associated with the kernel $K_A(x,y) = \frac{1}{\|x-y\|} \chi_{A}^c(x-y) = \kappa_A(x-y)$, for $x \not= y$, where $\chi_{A}^c$ denotes the indicator function of the cone in $\R^2$ of all the lines with directions in $A \subset S^1$. It is worth observing that the presence of a singularity in $\kappa$ along the diagonal (i.e., for $x=y$) depending only on the $-1$ power of $\|x-y\|$ precisely identifies $\mcR^* \mcR$ as a $\Psi$DO of order $-1$. Instead, the presence of discontinuities in $\kappa_A$ (also approaching the diagonal) depends on $x-y$ itself, and not only on its norm: this prevents its classification as a $\Psi$DO. As a consequence, the job 
of the correction $\PDONet$ should be replacing $\mcRA^* \mcRA$ with an operator whose Schwartz kernel is smooth.

To show that this is indeed the case, we first observe that: 
\[
\PDONet \colon L^2(\R^2)^Q \rightarrow L^2(\R^2)^Q, \quad [\PDONet w]_\iota = \sum_{\iota' \leq Q} \widetilde{\kappa}_{\iota,\iota'} \ast w_{\iota'},
\]
where $\widetilde{\kappa}_{\iota,\iota'} \in L^2(\R^2)$ are learnable functions. This is a direct consequence of ``lifting'' the original 
definition of the $\Psi$DONet architecture provided in \cite{bubba2021deep} to the continuous setting introduced in Section~\ref{ssec:ContFormulation}. The only negligible difference is that, in the discrete setting discussed in \cite{bubba2021deep}, upsampling and downsampling operators are employed to optimize the number of parameters needed to approximate the convolution operator via discrete convolutions. However, those operations are not needed in this continuous formulation, and here we avoid introducing them for simplicity.
We now show that it is  
possible to provide a similar representation for the operator $\mcRA^* \mcRA$ in the wavelet domain. Indeed, for $w,w' \in L^2(\R^2)^Q$, we have:
\[ \begin{aligned}
    \langle \mathscr{W} \mcRA^*\mcRA \mathscr{W}^* w, w' \rangle &= \int_{\R^2} \int_{\R^2} \kappa_A(x-y)\left(\sum_\iota \int_{\R^2}w_\iota(z)\psi_{\iota,z}(x)dz\right)\left(\sum_{\iota'} \int_{\R^2} w_{\iota'}(z)\psi_{\iota',z'}(x)dz'\right)dxdy
    \\
    & = \sum_{\iota,\iota'} \int_{\R^2}\int_{\R^2} K_A^{\iota,\iota'}(z,z') w_\iota(z)w_{\iota'}(z') dz dz',
    \end{aligned} 
\]
where each kernel, defined as
\[
K_A^{\iota,\iota'}(z,z') = \int_{\R^2} \int_{\R^2} \kappa_A(x-y) \psi_{\iota,z}(x)\psi_{\iota',z'}(y) dx dy
\]
can be interpreted as the wavelet subband of the wavelet transform associated with $K_A(x,y)$ with respect to both variables $x$ and $y$: $K_{A}^{\iota,\iota'} = \big[\mathscr{W}_y[\mathscr{W}_x K_A]_\iota\big]_{\iota'}$. Moreover, such kernels are convolutional, since (through the change of variables $\texttt{x} = x-z,\texttt{y} = y-z'$) we have:
\[
\begin{aligned}
    K_A^{\iota,\iota'}(z,z') &= \int_{\R^2} \int_{\R^2} \kappa_A(x-y) 2^{j+j'}\psi^{(t)}(2^{j}(x-z))\psi^{(t')}(2^{j'}(y-z')) dx dy \\
    &= \int_{\R^2} \int_{\R^2} \kappa_A(\texttt{x}-\texttt{y}+(z-z')) 2^{j+j'}\psi^{(t)}(2^{j}\texttt{x})\psi^{(t')}(2^{j'}\texttt{y}) d\texttt{x} d\texttt{y} \\[0.25em]
    &=: \kappa_A^{\iota,\iota'}(z-z'). 
\end{aligned}
\]
This allows to conclude that:
\begin{equation}\label{eq:PSOcorrection}
[\mathscr{W} \mcRA^*\mcRA \mathscr{W}^* w + \PDONet w]_\iota =  \sum_{\iota' \leq Q} (\kappa_A^{\iota,\iota'} + \widetilde{\kappa}_{\iota,\iota'})\ast w_{\iota'}.
\end{equation}
Therefore, the learned filters $\widetilde{\kappa}_{\iota,\iota'}$ succeed in modifying $\mcRA^*\mcRA$ at the core of its representation in the wavelet domain. Among the possible choices explored at training time, the learned filters might be such that, when added to $\kappa_A^{\iota,\iota'}$, they dampen the discontinuities introduced by $\chi_A^c$ in $\kappa_A$, thus allowing for 
smooth filters and, ultimately, a smoothed version of the operator $\mcRA^*\mcRA$.

\subsection{Optimizing $\Psi$DONet with special filters}\label{ssubsec:specialf}

In a discrete setting, the arguments discussed in Section \ref{ssubsec:kernels} have a clear interpretation. On the one hand, the learned convolutional kernels of each iteration of \eqref{eq:PsiDONet}, stored in $\zetab_k$, can be interpreted as a pixel-map image discretizing $\kappa^{\iota,\iota'}$. On the other hand, the action of the operator $\Wb \RbA^\top \RbA \Wb^\top$ can be viewed as a convolution with $S^2$ filters, which can be computed numerically following the procedure described in \cite[Section 5.2]{bubba2021deep}. In $\Psi$DONet, therefore, an artifact-suppressing correction is obtained by adding a learned convolutional filter to the ones describing $\Wb \RbA^\top \RbA \Wb^\top$ (cf.~Equation~\eqref{eq:PSOcorrection}). Such a filter, as a pixel-map image, is assumed to be supported in a (possibly small) square centred at the origin (cf.~Figure~\ref{fig:filters_geom}(a)). This choice reflects the $\Psi$DO nature of the normal operator $\mcR^*\mcR$, whose convolutional kernel $\kappa$ rapidly vanishes away from the origin, which suggests to concentrate the corrections closer to the centre and learn all of the parameters of the (possibly small) square filter.

Nevertheless, given the task of regularizing $\Wb \RbA^\top \RbA \Wb^\top$ by smoothing the singularities introduced by $\chi_{A}^c$ within the kernel $\kappa_A$ and its wavelet representations, we can resort to smartly crafted, masked filters to reduce the number of learnable parameters (cf. Figure~\ref{fig:filters_geom}). We therefore propose three variants of $\Psi$DONet, specifically suited for limited- and sparse-angle tomography.

\begin{itemize}
    \item \textbf{$\Psi$DONet-bow} (abbreviated in $\Psi_{\text{bow}}$). In the case of limited-angle tomography, $A = [-\Gamma,\Gamma]$, we first observe that the presence of $\chi_A^{c}$ in the filter $\kappa_A$ sets to zero all the elements outside the cone with angles from $-\Gamma$ to $\Gamma$. Therefore, it is unnecessary to learn any correction in those parts of the filters, and we can restrict the learnable parameters within the visible cone. When doing so, the learnable filters are no longer supported on a full square but on the intersection of the square with the cone, resembling the shape of a bowtie (cf.~Figure~\ref{fig:filters_geom}(b)). 
    \item \textbf{$\Psi$DONet-x} (abbreviated in $\Psi_{\text{x}}$). Still in the context of limited-angle tomography, we also realize that, in order to smooth out the jump of $\chi_A^{c}$ in $\kappa_A$ (or in its wavelet representation), it is not necessary to consider a correction supported in the whole visible cone, but it is sufficient to restrict it to a sufficiently large wedge close to the lines with directions in $\partial A = \pm \Gamma$. This implies that the learnable parameters form a cross shape, within a square, along the directions $\pm \Gamma$, with a small, fixed offset (cf.~Figure~\ref{fig:filters_geom}(c)). 
    \item \textbf{$\Psi$DONet-sparse} (abbreviated in $\Psi_{\text{spa}}$).  In the case of sparse-angle tomography, it is again meaningful to concentrate the smoothing correction close to the lines of directions in $\partial A$, which, in this case, is given by $\{\omega_1\pm \eta, \ldots \omega_N \pm \eta\}\approx \{\omega_1, \ldots, \omega_N\}$. The resulting learnable filters are supported within a square, along the lines of directions $\{\omega_1, \ldots, \omega_N\}$ and with a small offset (cf.~Figure~\ref{fig:filters_geom}(d)).  
\end{itemize}

\input{images/filters/filters_all}

\section{Numerical results}\label{sec:NumericalExperiments}
\input{numerical_results}

\section{Conclusions}\label{sec:Conclusions}
In this paper, we revisited $\Psi$DONet, a supervised learning approach based on unrolling the iterations of ISTA. Our main contribution is a novel microlocal interpretation for the theoretical analysis of $\Psi$DONet. In particular, leveraging the semi-discrete wavelet transform, we introduce a new continuous
formulation of $\Psi$DONet which, thanks to the very recent analysis provided in~\cite{andrade2022deep}, allows to
show that: (1) $\Psi$DONet might recover edges and discontinuities that are invisible from the incomplete sinogram, although we cannot predict their exact location; (2) $\Psi$DONet not only can dampen the strength of streak
artifacts, as initially shown in~\cite{bubba2021deep}, but can prevent them from appearing in the reconstruction. Compared to~\cite{andrade2022deep}, our theoretical analysis focuses on the analysis of the singular support rather than on the wavefront set: this is because our initial model uses wavelets, which are able to resolve only the singular support of a distribution. Further extending our analysis to the wavefront set would require the use of shearlets or curvelets, which in turn requires to consider a different primal-dual scheme for the solution of the variational formulation, and ultimately ``reinvent'' the $\Psi$DONet architecture. 

From a numerical standpoint, the microlocal interpretation motivates the introduction of three novel variants of the $\Psi$DONet architecture, based on smartly crafted filters whose structure is specifically inspired by the streak artifact singularities characterizing tomographic reconstructions from incomplete data. As a proof of concept, we tested the different versions of $\Psi$DONet on simulated data from limited- and sparse-angle geometry, generated from the ellipse data set, considering different visisble wedges and level of sparsity of the angular views. The results show that, despite the considerably lower number of (learnable) parameters, the quality of the reconstruction does not worsen compared to the case of square filters, and it is even improved in certain scenarios. 
Future extensions could address the limitations of finite unfolded architectures considering infinite depth models, deriving a theoretical characterization of the resulting reconstruction.

\section*{Acknowledgments}
The authors would like to thank the Isaac Newton Institute (INI) for Mathematical Sciences, Cambridge, for support and hospitality during the INI retreat in July 2024 where work on this paper was undertaken.
T.A.~Bubba and L.~Ratti acknowledge support by the European Union - NextGeneration EU through the Italian Ministry of University and Research as part of the PNRR – M4C2, Investment 1.3 (MUR Directorial Decree no. 341 of 03/15/2022), FAIR ``Future Partnership Artificial Intelligence Research'', Proposal Code PE00000013 - CUP J33C22002830006). 
A.~Sebastiani is supported by the project “PNRR - Missione 4 ``Istruzione e Ricerca'' - Componente C2 Investimento 1.1 ``Fondo per il Programma Nazionale di Ricerca e Progetti di Rilevante Interesse Nazionale (PRIN)'', ``Advanced optimization METhods for automated central veIn Sign detection in multiple sclerosis from magneTic resonAnce imaging (AMETISTA)'', project code: P2022J9SNP, MUR D.D. financing decree n.~1379 of 1st September 2023 (CUP E53D23017980001) funded by the European Commission under the NextGeneration EU programme.
A.~Sebastiani is partially supported by INdAM-GNCS (Project "MOdelli e MEtodi Numerici per il Trattamento delle Immagini (MOMENTI)", Code CUP\_E53C23001670001). 
T.A.~Bubba, L.~Ratti and A.~Sebastiani are members of INdAM-GNCS.

\bibliographystyle{plain}
\bibliography{References.bib}

\end{document}

%% file: images/filters/filters_all.tex
\begin{figure}
	\centering
	\subfloat[$\Psi_{\text{o}}$ filters]{
		\centering
		\begin{tikzpicture}[scale=0.3]
		\draw[help lines] (-5,-5) grid (6,6);
            \fill[cyan!80] (-3,-3) -- (-3, 4) -- (4,4) -- (4,-3);
	    \end{tikzpicture}
	}
        \hfill
	\subfloat[$\Psi_{\text{bow}}$ filters]{
		\centering
		\begin{tikzpicture}[scale=0.3]
		\draw[help lines] (-5,-5) grid (6,6);
            \fill[cyan!80] (0.5,0.5) -- (-3,-3.75) -- (-3,4.75);
            \fill[cyan!80] (0.5,0.5) -- (4,4.75) -- (4,-3.75);
            \draw[cyan!50] (-4, -5) -- (5, 6);
            \draw[cyan!50] (-4, 6) -- (5, -5);
		\end{tikzpicture}
	}
        \hfill
	\subfloat[$\Psi_{\text{x}}$ filters]{
		\centering
		\begin{tikzpicture}[scale=0.3]
		\draw[help lines] (-5,-5) grid (6,6);	
            \draw[cyan!50] (-4, -5) -- (5, 6);
            \draw[cyan!50] (-4, 6) -- (5, -5);
            \draw[cyan!80, line width=0.85mm] (-3, -3.75) -- (4, 4.75);
            \draw[cyan!80, line width=0.85mm] (-3, 4.75) -- (4, -3.75);
	    \end{tikzpicture}
	}
        \hfill
        \subfloat[$\Psi_{\text{spa}}$ filters]{
		\centering
		\begin{tikzpicture}[scale=0.3]
		\draw[help lines] (-5,-5) grid (6,6);	
            \draw[cyan!50] (-2.7, -5) -- (3.7, 6);
            \draw[cyan!50] (-2.7, 6) -- (3.7, -5);
            \draw[cyan!50] (-5, -2.7) -- (6, 3.7);
            \draw[cyan!50] (6, -2.7) -- (-5, 3.7);
            
            \draw[cyan!50] (0.5, 6) -- (0.5, -5);
            \draw[cyan!50] (6, 0.5) -- (-5, 0.5);
            
            \draw[cyan!80, line width=0.85mm] (-1.5, -3) -- (2.5, 4);
            \draw[cyan!80, line width=0.85mm] (-1.5, 4) -- (2.5, -3);
            \draw[cyan!80, line width=0.85mm] (-3, -1.5) -- (4, 2.5);
            \draw[cyan!80, line width=0.85mm] (4, -1.5) -- (-3, 2.5);
            \draw[cyan!80, line width=0.85mm] (-3, 0.5) -- ( 4, 0.5);
            \draw[cyan!80, line width=0.85mm] (0.5,-3) -- (0.5, 4);
	    \end{tikzpicture}
	}
    \caption{Different filter geometries for $\Psi$DONet. (a) $\Psi$DONet filters initially proposed in~\cite{bubba2021deep}. (b) $\Psi$DONet-bow filters for limited-angle tomography. (c) $\Psi$DONet-x filters for limited-angle tomography. (d) $\Psi$DONet-sparse filters for sparse-angle tomography. } 
\label{fig:filters_geom}
\end{figure}
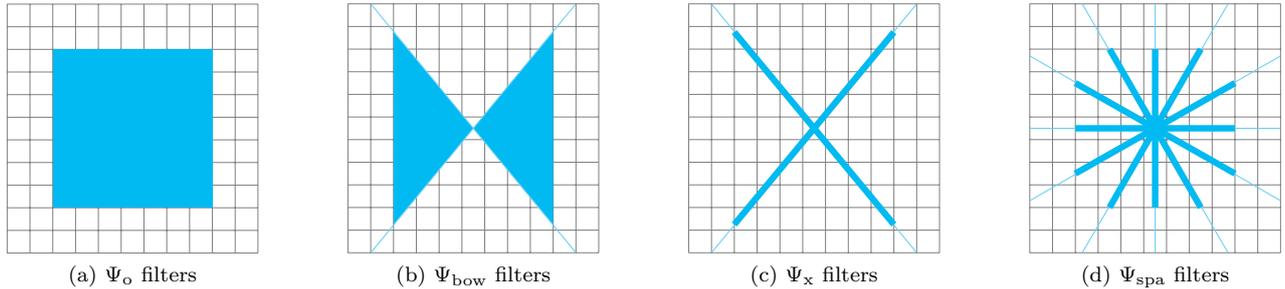

%% file: numerical_results.tex
\newcommand{\uu}{{\ub}}
\newcommand{\uupsi}{{\ub_{\Psi_{\text{o}}}}}

\newcommand{\uubow}[1]{{\ub_{\Psi_{\text{bow}}^{#1}}}}
\newcommand{\uux}[1]{{\ub_{\Psi_{\text{x}}^{#1}}}}
\newcommand{\uuspa}[1]{{\ub_{\Psi_{\text{spa}}^{#1}}}}
\newcommand*{\discreteop}{\RbA}

In this section, we evaluate the performance of the proposed variants $\Psi_{\text{bow}}$, $\Psi_{\text{x}}$ and $\Psi_{\text{spa}}$ of $\Psi$DONet by comparing their performances with that of the original $\Psi$DONet (later on simply denoted by $\Psi_{\text{o}}$) introduced in~\cite{bubba2021deep}. This mostly serves as a proof of concept that reducing the support of the filters, thanks to the insights coming from microlocal analysis, does not hinder the quality of the reconstruction, despite the lower number of learnable parameters.

We first describe the implementation\footnote{The code to reproduce all results in this paper will be made available upon acceptance of this manuscript.} of $\Psi$DONet. The code released in~\cite{bubba2021deep} has been updated in several parts, to introduce novel efficient libraries for the computation of the involved operators, to keep up with the latest release of some packages and to implement the novel filters. Most importantly, the implementation is now in PyTorch. For the sake of completeness, we point out that the repository\footnote{\url{https://github.com/megalinier/PsiDONet}} associated with~\cite{bubba2021deep} also offers a PyTorch version of the Tensorflow codes which uses the ODL library \cite{adler2017operator}. Instead, the current implementation relies on Tomosipo \cite{hendriksen2021tomosipo} for the discretization of the subsampled Radon transform $\discreteop$. In our preliminary experiments this major change significantly reduced the time required for the training. The CT acquisition geometry is specified
by means of the library LION\footnote{\url{https://github.com/CambridgeCIA/LION}} and we further normalize both $\discreteop$ and $\discreteop^\top$ to ensure they have unitary norm. As in the original implementation, we use PyWavelets~\cite{lee2019pywavelets} to construct the wavelet transform. In our experiments, we consider the Haar wavelet basis \cite{mallat1999wavelet} setting the maximum and minimum scales to $J=7$ and $J_0=4$, thus having a wavelet decomposition with $Q=10$ subbands. 

Concerning the architecture, we fix the number of unrolled blocks to 10. The convolutional filters of the variants $\Psi_{\text{bow}}$, $\Psi_{\text{x}}$ and $\Psi_{\text{spa}}$ are obtained by multiplying the weights of the full square convolutional filter with a boolean mask. The construction process of these boolean masks can be described by referring to Figure \ref{fig:filters_geom}. Assuming that the dimensions of the convolutional filter are odd, the pixel at the center is considered as the origin, and the pixels intersecting the lines, whose orientation is contained in $A$ (or $\partial A$), as described in Section \ref{ssubsec:specialf}, are included in the mask. This discretization may introduce some inaccuracies depending on the size of the filters. For this reason, in our experiments we modify these masks considering stripes having a half-width of $q=1,2,3$ pixels in place of lines of 1-pixel width. Notice that the use of a mask ensures that the learnable parameters are only those in the support specified by the filter geometry, however, it does not necessarily come with a considerable speed-up for training times, which is still dominated by the size of the convolutional filters.

For training and testing, we use the same ellipse datasets released in \cite{bubba2021deep}. More in detail, they consist of synthetic images of size $128\times 128$ containing randomly generated ellipses with different locations, sizes, numbers, and intensity gradients. The training and testing datasets contain 10000 and 500 images, respectively. The simulated measurements available in the datasets are already provided including additive Gaussian noise. All the models are trained using the Adam optimizer setting the learning rate at $10^{-3}$. The batch size is set to $25$ and the number of training epochs is set at to $15$. The final reconstructions are assessed in terms of peak signal-to-noise ratio (PSNR) and structured similarity index (SSIM). All the experiments have been performed on a workstation equipped with an Intel i9-12900K processor and an NVIDIA RTX A5000 GPU.

\subsection{Limited-angle tomography}
In this section we focus on the limited-angle CT problem considering two different settings, corresponding to missing wedges of $60^{\circ}$ and $120^{\circ}$ in parallel-beam geometry. In both the experiments, we set the size of the convolutional filters to $33 \times 33$. 
For the sake of readability, we use the following abbreviations to denote the results of the different reconstruction models:
\begin{itemize}
    \item $\uupsi$: standard $\Psi$DONet with full convolutional filters $\Psi_{\text{o}}$.
    \item $\uubow{0}$: $\Psi$DONet with $\Psi_{\text{bow}}$  filters with 1-pixel width lines in the boolean mask. 
    \item $\uubow{q}$: $\Psi$DONet with $\Psi_{\text{bow}}$  filters with stripes of half-width $q$ in the boolean mask.
    \item $\uux{0}$: $\Psi$DONet with $\Psi_{\text{x}}$ filters with 1-pixel width lines in the boolean mask.
    \item $\uux{q}$: $\Psi$DONet with $\Psi_{\text{x}}$ filters with stripes of half-width $q$ in the boolean mask.
\end{itemize}

For the first experiment, we consider the set of angles $A=[-\pi/3,\pi/3]$ corresponding to a missing wedge of $60^{\circ}$. We report in Table \ref{tab:LA30} the average PSNR and SSIM computed on the test set images. 
All the values are really close to each other, confirming the benefit of learning the correction only on the support specified by $\chi_A^{c}$. However, we observe that the 1-pixel width variants, $\Psi_{\text{bow}}^0$ and $\Psi_{\text{x}}^0$ results in slightly worse values, confirming the limitation caused by the line discretization. 
In contrast, the use of the ``stripe model'' shows an improvement, for both the $\Psi_{\text{bow}}$ and $\Psi_{\text{x}}$ setting. Finally, the metrics for $\Psi_{\text{x}}$ show that expanding the support over a certain threshold introduces numerical instabilities. Overall, these results confirm that restricting the filters support does not significantly compromise the quality of the reconstruction.

\begin{table}[h]
    \centering
    \begin{tabular}{|c|cccccc|}
    \cline{2-7}\multicolumn{1}{c|}{}
        & $\uupsi$ & $\uubow{0}$ & $\uubow{3}$ & $\uux{0}$ & $\uux{1}$ & $\uux{3}$ \\
    \hline
PSNR & 29.01 & 28.89 & 29.05 & 28.40 & 29.19 & 29.13 \\
SSIM & 0.872 & 0.868 & 0.872 & 0.852 & 0.876 & 0.870 \\
    \hline
    \end{tabular}
    \caption{Average PSNR and SSIM on the test set for limited-angle with $A = [-\pi/3,\pi/3]$. Sample reconstructions can be found in Figure~\ref{fig:LA30_mixed}.}
    \label{tab:LA30}
\end{table}
In Figure \ref{fig:LA30_mixed} we show the reconstructions of two images from the test set having different distributions of the ellipses. The reconstructions are very similar and we highlight the differences considering two close-ups. In general, it can be observed that the $\Psi_{x}$ variants provide slightly better reconstructions in terms of smoothness and noise suppression. We can see that small objects appear with neat boundaries, without the presence of streak artifacts.

\input{images/LA30/mixed}

Next, we consider the set of angles $A=[-\pi/6,\pi/6]$, corresponding to a missing wedge of $120^{\circ}$. We report the average PSNR and SSIM in Table \ref{tab:LA30}. Despite having reduced the amount of available data, we can draw similar considerations to those of the case $A=[-\pi/3,\pi/3]$. The average values of the quality metrics are all close to each other, with only the 1-pixel width variant $\Psi_{\text{x}}^0$ resulting in slightly smaller values. 

\begin{table}[h]
    \centering
    \begin{tabular}{|c|cccccc|}
    \cline{2-7}
        \multicolumn{1}{c|}{} & $\uupsi$ & $\uubow{0}$ & $\uubow{3}$ & $\uux{0}$ & $\uux{1}$ & $\uux{3}$\\
    \hline
PSNR & 25.33 & 25.51 & 25.57 & 25.02 & 25.65 & 25.65\\
SSIM & 0.754 & 0.755 & 0.760 & 0.738 & 0.771 & 0.771 \\
    \hline
    \end{tabular}
    \caption{Average PSNR and SSIM on the test set for limited-angle with $A=[-\pi/6,\pi/6]$. Sample reconstructions can be found in Figure~\ref{fig:LA60_mixed}.}
    \label{tab:LA60}
\end{table}

In Figure \ref{fig:LA60_mixed} we show the reconstructions of two images from the test set, having different distributions of the ellipses. From the close-ups in the top row, it is evident that the smaller values of quality metrics with $\Psi_{\text{o}}$ filters is likely due to the suppressed tiny and low-contrast regions, which deform the original shape boundaries. Instead, the restricted versions of the filters $\Psi_{\text{bow}}$ and $\Psi_{\text{x}}$ are able to shade the presence even of the smaller features. From the bottom row of Figure \ref{fig:LA60_mixed} we can observe that all the methods deliver almost indistinguishable reconstructions, to the point that the small ellipse at the center of the image is not visible in any of the reconstructions. This suggests that we cannot control nor predict which new singularities (not present in the data) $\Psi$-DONet is able to introduce in the final reconstruction. 

\input{images/LA60/mixed}

\subsection{Sparse-angle tomography}
In this section, we test $\Psi$DONet on sparse-angle CT problems, considering 12 and 6 uniformly spaced angles in the set $[-\pi/2,\pi/2]$ in parallel-beam geometry, which corresponds to taking an angle every $30^{\circ}$ and  $60^{\circ}$, respectively. We point out that these settings are very challenging, since the information available in the measurements is very scarce.
In this case, building the support for the mask needs to be handled with care.
In a preliminary test, we observed that the slices defining the $\Psi_{\text{spa}}$ boolean masks cover the entire support when considering filters of size $33\times 33$, as in the experiments of the previous section.
For this reason, motivated by the experimental analysis on the approximation error of the convolution representation of $\Wb \RbA^\top \RbA \Wb^\top$ depending on the filters size, reported in \cite[Section 5.2 and Figure 7]{bubba2021deep}, we consider convolutional filters of size $65\times 65$ for the experiments in this section, to be able to better examine the potential benefit of reducing the filter support.

We report in Table \ref{tab:sparse_angles} the average PSNR and SSIM computed on the test set for both the cases, denoting by $\uuspa{1}$ the reconstruction obtained using $\Psi$DONet with the $\Psi_{\text{spa}}$ filters whose boolean masks also include the stripes having half-width of $q=1$ pixel. It is evident that the use of filters with smaller supports does not worsen the performance of the method. In addition, there is a small improvement in the metrics or, equivalently, a reduction of the instabilities due to the parametrization of the ``stripe'' version of the filters.

\begin{table}[h]
    \centering
    \begin{tabular}{|c|cc|cc|}
        \cline{2-5}
         \multicolumn{1}{c|}{} & \multicolumn{2}{|c|}{12 angles} & \multicolumn{2}{|c|}{6 angles}\\
        \cline{2-5}
          \multicolumn{1}{c|}{} & $\uupsi$ & $\uuspa{1}$ & $\uupsi$ & $\uuspa{1}$ \\
        \hline
         PSNR & 28.27 & 28.30 & 25.43 & 25.59 \\
         SSIM &  0.833 & 0.834 & 0.735 & 0.737 \\
        \hline
    \end{tabular}
    \caption{Average PSNR and SSIM on the test set sparse angle, considering 12 and 6 angles. Sample reconstructions can be found in Figure~\ref{fig:sparse15_mixed} and Figure \ref{fig:sparse30_mixed}, respectively.}
    \label{tab:sparse_angles}
\end{table}

We show in Figure \ref{fig:sparse15_mixed} and Figure \ref{fig:sparse30_mixed} the reconstructions for two sets of images from the test set using 12 and 6 angles, respectively. First of all, we notice how despite the very undersampled data, the overall reconstruction does not show evident regions where the information is missing nor clear streak artifacts. More in detail, in Figure \ref{fig:sparse15_mixed} we observed that the use of $\Psi_{\text{spa}}$ filters slightly improves the smoothness of constant regions and the edges in overlapping ellipses. Concerning the reconstruction shown in Figure \ref{fig:sparse30_mixed}, 
we can sometimes spot the appearance of some (partial) streak artifacts, possibly due to the highly reduced amount of measured data. Despite that, the reconstructions obtained with $\Psi_{\text{spa}}$ filters show a reduced magnitude of these distortions. Moreover, as in the other setting, the $\Psi_{\text{spa}}$ filters seem to better preserve the boundaries in overlapping regions.

\input{images/sparse15/mixed}
\input{images/sparse30/mixed}

%% file: images/LA30/mixed.tex
\begin{figure}
	\centering
	\subfloat[GT]{
        \begin{tabular}{c}
	   \scalebox{0.8}{
    	\begin{tikzpicture}
    	\begin{scope}[spy using outlines={rectangle,blue,magnification=2,size=1.5cm}]
    	\node [name=c]{	\includegraphics[height=3cm]{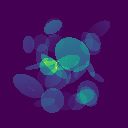}};
    	\spy on (0.8,0.4) in node [name=c1]  at (0.75,-2.25);
    	\spy on (-0.5,0.2) in node [name=c1]  at (-0.75,-2.25);
    	\end{scope}
    	\end{tikzpicture}}\\
            \scalebox{0.8}{
    	\begin{tikzpicture}
    	\begin{scope}[spy using outlines={rectangle,blue,magnification=2,size=1.5cm}]
    	\node [name=c]{	\includegraphics[height=3cm]{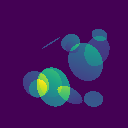}};
    	\spy on (0.8,0.4) in node [name=c1]  at (0.75,-2.25);
    	\spy on (-0.5,0.2) in node [name=c1]  at (-0.75,-2.25);
    	\end{scope}
    	\end{tikzpicture}}
        \end{tabular}
        }
	\subfloat[$\uupsi$]{
        \begin{tabular}{c}
	\scalebox{0.8}{
	\begin{tikzpicture}
	\begin{scope}[spy using outlines={rectangle,blue,magnification=2,size=1.5cm}]
	\node [name=c]{	\includegraphics[height=3cm]{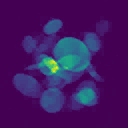}};
	\spy on (0.8,0.4) in node [name=c1]  at (0.75,-2.25);
	\spy on (-0.5,0.2) in node [name=c1]  at (-0.75,-2.25);
	\end{scope}
	\end{tikzpicture}}\\
	\scalebox{0.8}{
	\begin{tikzpicture}
	\begin{scope}[spy using outlines={rectangle,blue,magnification=2,size=1.5cm}]
	\node [name=c]{	\includegraphics[height=3cm]{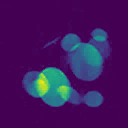}};
	\spy on (0.8,0.4) in node [name=c1]  at (0.75,-2.25);
	\spy on (-0.5,0.2) in node [name=c1]  at (-0.75,-2.25);
	\end{scope}
	\end{tikzpicture}}
        \end{tabular}
        }
	\subfloat[$\uubow{0}$]{
        \begin{tabular}{c}
	\scalebox{0.8}{
	\begin{tikzpicture}
	\begin{scope}[spy using outlines={rectangle,blue,magnification=2,size=1.5cm}]
	\node [name=c]{	\includegraphics[height=3cm]{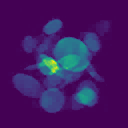}};
	\spy on (0.8,0.4) in node [name=c1]  at (0.75,-2.25);
	\spy on (-0.5,0.2) in node [name=c1]  at (-0.75,-2.25);
	\end{scope}
	\end{tikzpicture}}\\
        \scalebox{0.8}{
	\begin{tikzpicture}
	\begin{scope}[spy using outlines={rectangle,blue,magnification=2,size=1.5cm}]
	\node [name=c]{	\includegraphics[height=3cm]{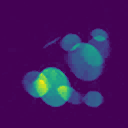}};
	\spy on (0.8,0.4) in node [name=c1]  at (0.75,-2.25);
	\spy on (-0.5,0.2) in node [name=c1]  at (-0.75,-2.25);
	\end{scope}
	\end{tikzpicture}}
        \end{tabular}
        }
	\subfloat[$\uux{0}$]{
        \begin{tabular}{c}
	\scalebox{0.8}{
	\begin{tikzpicture}
	\begin{scope}[spy using outlines={rectangle,blue,magnification=2,size=1.5cm}]
	\node [name=c]{	\includegraphics[height=3cm]{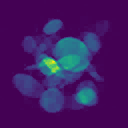}};
	\spy on (0.8,0.4) in node [name=c1]  at (0.75,-2.25);
	\spy on (-0.5,0.2) in node [name=c1]  at (-0.75,-2.25);
	\end{scope}
	\end{tikzpicture}}\\
        \scalebox{0.8}{
	\begin{tikzpicture}
	\begin{scope}[spy using outlines={rectangle,blue,magnification=2,size=1.5cm}]
	\node [name=c]{	\includegraphics[height=3cm]{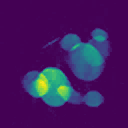}};
	\spy on (0.8,0.4) in node [name=c1]  at (0.75,-2.25);
	\spy on (-0.5,0.2) in node [name=c1]  at (-0.75,-2.25);
	\end{scope}
	\end{tikzpicture}}
        \end{tabular}}
	\subfloat[$\uux{1}$]{
        \begin{tabular}{c}
	\scalebox{0.8}{
	\begin{tikzpicture}
	\begin{scope}[spy using outlines={rectangle,blue,magnification=2,size=1.5cm}]
	\node [name=c]{	\includegraphics[height=3cm]{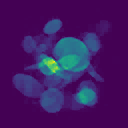}};
	\spy on (0.8,0.4) in node [name=c1]  at (0.75,-2.25);
	\spy on (-0.5,0.2) in node [name=c1]  at (-0.75,-2.25);
	\end{scope}
	\end{tikzpicture}}\\
        \scalebox{0.8}{
	\begin{tikzpicture}
	\begin{scope}[spy using outlines={rectangle,blue,magnification=2,size=1.5cm}]
	\node [name=c]{	\includegraphics[height=3cm]{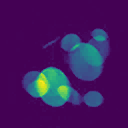}};
	\spy on (0.8,0.4) in node [name=c1]  at (0.75,-2.25);
	\spy on (-0.5,0.2) in node [name=c1]  at (-0.75,-2.25);
	\end{scope}
	\end{tikzpicture}}
        \end{tabular}}	
    \caption{Visualization of the reconstructions of two images  from the test set with $A=[-\pi/3,\pi/3]$.}
    \label{fig:LA30_mixed}
\end{figure}

%% file: images/LA60/mixed.tex
\begin{figure}
	\centering
	\subfloat[GT]{
        \begin{tabular}{c}
	   \scalebox{0.8}{
    	\begin{tikzpicture}
	\begin{scope}[spy using outlines={rectangle,blue,magnification=2,size=1.5cm}]
	\node [name=c]{	\includegraphics[height=3cm]{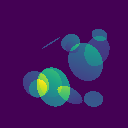}};
	\spy on (0.5,-0.75) in node [name=c1]  at (0.75,-2.25);
	\spy on (-0.2,0.35)  in node [name=c1]  at (-0.75,-2.25);
	\end{scope}
	\end{tikzpicture}}\\
            \scalebox{0.8}{
    	\begin{tikzpicture}
	\begin{scope}[spy using outlines={rectangle,blue,magnification=2,size=1.5cm}]
	\node [name=c]{	\includegraphics[height=3cm]{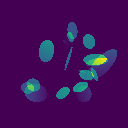}};
	\spy on (0.85,0.35) in node [name=c1]  at (0.75,-2.25);
	\spy on (-0.6,-0.6)  in node [name=c1]  at (-0.75,-2.25);
	\end{scope}
	\end{tikzpicture}}
        \end{tabular}
        }
	\subfloat[$\uupsi$]{
        \begin{tabular}{c}
	\scalebox{0.8}{
	\begin{tikzpicture}
	\begin{scope}[spy using outlines={rectangle,blue,magnification=2,size=1.5cm}]
	\node [name=c]{	\includegraphics[height=3cm]{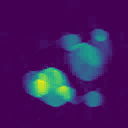}};
	\spy on (0.5,-0.75) in node [name=c1]  at (0.75,-2.25);
	\spy on (-0.2,0.35)  in node [name=c1]  at (-0.75,-2.25);
	\end{scope}
	\end{tikzpicture}}\\
	\scalebox{0.8}{
	\begin{tikzpicture}
	\begin{scope}[spy using outlines={rectangle,blue,magnification=2,size=1.5cm}]
	\node [name=c]{	\includegraphics[height=3cm]{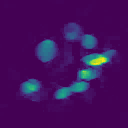}};
	\spy on (0.85,0.35) in node [name=c1]  at (0.75,-2.25);
	\spy on (-0.6,-0.6)  in node [name=c1]  at (-0.75,-2.25);
	\end{scope}
	\end{tikzpicture}}
        \end{tabular}
        }
	\subfloat[$\uubow{0}$]{
        \begin{tabular}{c}
	\scalebox{0.8}{
	\begin{tikzpicture}
	\begin{scope}[spy using outlines={rectangle,blue,magnification=2,size=1.5cm}]
	\node [name=c]{	\includegraphics[height=3cm]{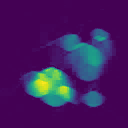}};
	\spy on (0.5,-0.75) in node [name=c1]  at (0.75,-2.25);
	\spy on (-0.2,0.35)  in node [name=c1]  at (-0.75,-2.25);
	\end{scope}
	\end{tikzpicture}}\\
        \scalebox{0.8}{
	\begin{tikzpicture}
	\begin{scope}[spy using outlines={rectangle,blue,magnification=2,size=1.5cm}]
	\node [name=c]{	\includegraphics[height=3cm]{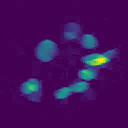}};
	\spy on (0.85,0.35) in node [name=c1]  at (0.75,-2.25);
	\spy on (-0.6,-0.6)  in node [name=c1]  at (-0.75,-2.25);
	\end{scope}
	\end{tikzpicture}}
        \end{tabular}
        }
	\subfloat[$\uux{0}$]{
        \begin{tabular}{c}
	\scalebox{0.8}{
	\begin{tikzpicture}
	\begin{scope}[spy using outlines={rectangle,blue,magnification=2,size=1.5cm}]
	\node [name=c]{	\includegraphics[height=3cm]{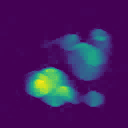}};
	\spy on (0.5,-0.75) in node [name=c1]  at (0.75,-2.25);
	\spy on (-0.2,0.35)  in node [name=c1]  at (-0.75,-2.25);
	\end{scope}
	\end{tikzpicture}}\\
        \scalebox{0.8}{
	\begin{tikzpicture}
	\begin{scope}[spy using outlines={rectangle,blue,magnification=2,size=1.5cm}]
	\node [name=c]{	\includegraphics[height=3cm]{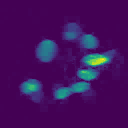}};
	\spy on (0.85,0.35) in node [name=c1]  at (0.75,-2.25);
	\spy on (-0.6,-0.6)  in node [name=c1]  at (-0.75,-2.25);
	\end{scope}
	\end{tikzpicture}}
        \end{tabular}}
	\subfloat[$\uux{1}$]{
        \begin{tabular}{c}
	\scalebox{0.8}{
	\begin{tikzpicture}
	\begin{scope}[spy using outlines={rectangle,blue,magnification=2,size=1.5cm}]
	\node [name=c]{	\includegraphics[height=3cm]{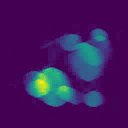}};
	\spy on (0.5,-0.75) in node [name=c1]  at (0.75,-2.25);
	\spy on (-0.2,0.35)  in node [name=c1]  at (-0.75,-2.25);
	\end{scope}
	\end{tikzpicture}}\\
        \scalebox{0.8}{
	\begin{tikzpicture}
	\begin{scope}[spy using outlines={rectangle,blue,magnification=2,size=1.5cm}]
	\node [name=c]{	\includegraphics[height=3cm]{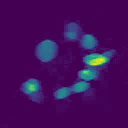}};
	\spy on (0.85,0.35) in node [name=c1]  at (0.75,-2.25);
	\spy on (-0.6,-0.6)  in node [name=c1]  at (-0.75,-2.25);
	\end{scope}
	\end{tikzpicture}}
        \end{tabular}}	
    \caption{Visualization of the reconstructions of two images from the test set with $A=[-\pi/6,\pi/6]$.}
    \label{fig:LA60_mixed}
\end{figure}

%% file: images/sparse15/mixed.tex
\begin{figure}
	\centering
	\subfloat[GT]{
        \begin{tabular}{c}
	   \scalebox{0.8}{
    	\begin{tikzpicture}
	\begin{scope}[spy using outlines={rectangle,blue,magnification=2,size=1.5cm}]
	\node [name=c]{	\includegraphics[height=3cm]{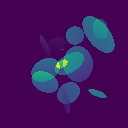}};
	\spy on (0.85,0.65) in node [name=c1]  at (0.75,-2.25);
	\spy on (-0.6,-0.6)  in node [name=c1]  at (-0.75,-2.25);
	\end{scope}
	\end{tikzpicture}}\\
            \scalebox{0.8}{
    	\begin{tikzpicture}
	\begin{scope}[spy using outlines={rectangle,blue,magnification=2,size=1.5cm}]
	\node [name=c]{	\includegraphics[height=3cm]{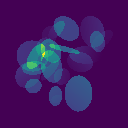}};
	\spy on (0.85,0.65) in node [name=c1]  at (0.75,-2.25);
	\spy on (-0.6,-0.6)  in node [name=c1]  at (-0.75,-2.25);
	\end{scope}
	\end{tikzpicture}}
        \end{tabular}
        }
	\subfloat[$\uupsi$]{
        \begin{tabular}{c}
	\scalebox{0.8}{
	\begin{tikzpicture}
	\begin{scope}[spy using outlines={rectangle,blue,magnification=2,size=1.5cm}]
	\node [name=c]{	\includegraphics[height=3cm]{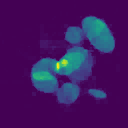}};
	\spy on (0.85,0.65) in node [name=c1]  at (0.75,-2.25);
	\spy on (-0.6,-0.6)  in node [name=c1]  at (-0.75,-2.25);
	\end{scope}
	\end{tikzpicture}}\\
	\scalebox{0.8}{
	\begin{tikzpicture}
	\begin{scope}[spy using outlines={rectangle,blue,magnification=2,size=1.5cm}]
	\node [name=c]{	\includegraphics[height=3cm]{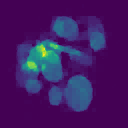}};
	\spy on (0.85,0.65) in node [name=c1]  at (0.75,-2.25);
	\spy on (-0.6,-0.6)  in node [name=c1]  at (-0.75,-2.25);
	\end{scope}
	\end{tikzpicture}}
        \end{tabular}
        }
	\subfloat[$\uuspa{1}$]{
        \begin{tabular}{c}
	\scalebox{0.8}{
	\begin{tikzpicture}
	\begin{scope}[spy using outlines={rectangle,blue,magnification=2,size=1.5cm}]
	\node [name=c]{	\includegraphics[height=3cm]{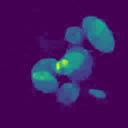}};
	\spy on (0.85,0.65) in node [name=c1]  at (0.75,-2.25);
	\spy on (-0.6,-0.6)  in node [name=c1]  at (-0.75,-2.25);
	\end{scope}
	\end{tikzpicture}}\\
        \scalebox{0.8}{
	\begin{tikzpicture}
	\begin{scope}[spy using outlines={rectangle,blue,magnification=2,size=1.5cm}]
	\node [name=c]{	\includegraphics[height=3cm]{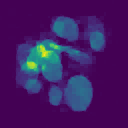}};
	\spy on (0.85,0.65) in node [name=c1]  at (0.75,-2.25);
	\spy on (-0.6,-0.6)  in node [name=c1]  at (-0.75,-2.25);
	\end{scope}
	\end{tikzpicture}}
        \end{tabular}
        }
    \caption{Visualization of the reconstructions of two images with 12 angles.}
    \label{fig:sparse15_mixed}
\end{figure}

%% file: images/sparse30/mixed.tex
\begin{figure}
	\centering
	\subfloat[GT]{
        \begin{tabular}{c}
	   \scalebox{0.8}{
    	\begin{tikzpicture}
	\begin{scope}[spy using outlines={rectangle,blue,magnification=2,size=1.5cm}]
	\node [name=c]{	\includegraphics[height=3cm]{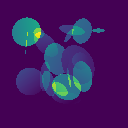}};
	\spy on (0.55,-0.55) in node [name=c1]  at (0.75,-2.25);
	\spy on (-0.6,0.3)  in node [name=c1]  at (-0.75,-2.25);
	\end{scope}
	\end{tikzpicture}}\\
            \scalebox{0.8}{
    	\begin{tikzpicture}
	\begin{scope}[spy using outlines={rectangle,blue,magnification=2,size=1.5cm}]
	\node [name=c]{	\includegraphics[height=3cm]{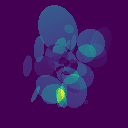}};
	\spy on (0.75,0.55) in node [name=c1]  at (0.75,-2.25);
	\spy on (-0.6,-0.6)  in node [name=c1]  at (-0.75,-2.25);
	\end{scope}
	\end{tikzpicture}}
        \end{tabular}
        }
	\subfloat[$\uupsi$]{
        \begin{tabular}{c}
	\scalebox{0.8}{
	\begin{tikzpicture}
	\begin{scope}[spy using outlines={rectangle,blue,magnification=2,size=1.5cm}]
	\node [name=c]{	\includegraphics[height=3cm]{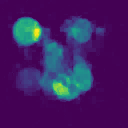}};
	\spy on (0.55,-0.55) in node [name=c1]  at (0.75,-2.25);
	\spy on (-0.6,0.3)  in node [name=c1]  at (-0.75,-2.25);
	\end{scope}
	\end{tikzpicture}}\\
	\scalebox{0.8}{
	\begin{tikzpicture}
	\begin{scope}[spy using outlines={rectangle,blue,magnification=2,size=1.5cm}]
	\node [name=c]{	\includegraphics[height=3cm]{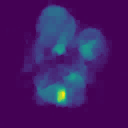}};
	\spy on (0.75,0.55) in node [name=c1]  at (0.75,-2.25);
	\spy on (-0.6,-0.6)  in node [name=c1]  at (-0.75,-2.25);
	\end{scope}
	\end{tikzpicture}}
        \end{tabular}
        }
	\subfloat[$\uuspa{1}$]{
        \begin{tabular}{c}
	\scalebox{0.8}{
	\begin{tikzpicture}
	\begin{scope}[spy using outlines={rectangle,blue,magnification=2,size=1.5cm}]
	\node [name=c]{	\includegraphics[height=3cm]{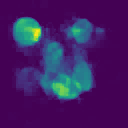}};
	\spy on (0.55,-0.55) in node [name=c1]  at (0.75,-2.25);
	\spy on (-0.6,0.3)  in node [name=c1]  at (-0.75,-2.25);
	\end{scope}
	\end{tikzpicture}}\\
        \scalebox{0.8}{
	\begin{tikzpicture}
	\begin{scope}[spy using outlines={rectangle,blue,magnification=2,size=1.5cm}]
	\node [name=c]{	\includegraphics[height=3cm]{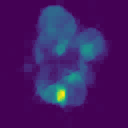}};
	\spy on (0.75,0.55) in node [name=c1]  at (0.75,-2.25);
	\spy on (-0.6,-0.6)  in node [name=c1]  at (-0.75,-2.25);
	\end{scope}
	\end{tikzpicture}}
        \end{tabular}
        }
    \caption{Visualization of the reconstructions of two images with 6 angles.}
    \label{fig:sparse30_mixed}
\end{figure}